\documentclass{article}
\usepackage{bm}
\usepackage{amsfonts, amsmath, amssymb}
\usepackage{lineno}
\usepackage{placeins} 
\usepackage[margin=1in]{geometry}
\usepackage[affil-it]{authblk} 
\usepackage{stmaryrd} 
\usepackage{setspace} 
\usepackage{floatrow} 
\usepackage{subfig} 
\usepackage[export]{adjustbox} 
\usepackage{sidecap} 
\usepackage[margin=1cm]{caption} 
\usepackage{graphicx}
\usepackage[labelsep=period]{caption} 
\usepackage{cite} 
\usepackage{times} 
\usepackage{mathptmx} 
\usepackage{booktabs,multirow}
\usepackage{microtype} 
\usepackage[colorlinks=true,linkcolor=blue,urlcolor=black,bookmarksopen=true]{hyperref}
\usepackage{bookmark}
\usepackage{tabularray}
\usepackage{longtable}
\usepackage{array}
\usepackage[title]{appendix}

\newcommand\keywords[1]%
  {\begin{flushleft}
   \let\and\\%
   \textbf{Keywords:}
   #1
   \end{flushleft}%
  }

\begin{document}

\title{Time integration schemes based on neural networks for solving partial differential equations on coarse grids}

\author[1,2]{Xinxin Yan}
\author[1,2]{Zhideng Zhou}
\author[3]{Xiaohan Cheng}
\author[1,2,*]{Xiaolei Yang}

\renewcommand\Affilfont{\fontsize{9}{10.8}\itshape}
\affil[1]{The State Key Laboratory of Nonlinear Mechanics, Institute of Mechanics, Chinese Academy of Sciences, \protect\\Beijing 100190, China.}
\affil[2]{School of Engineering Sciences, University of Chinese Academy of Sciences, Beijing 100049, China.}
\affil[3]{School of Science,Chang 'an University, Xi 'an, 710064,China}
\affil[*]{Corresponding author at xyang@imech.ac.cn}
\maketitle
\abstract
The accuracy of solving partial differential equations (PDEs) on coarse grids is greatly affected by the choice of discretization schemes. In this work, we propose to learn time integration schemes based on neural networks which satisfy three distinct sets of mathematical constraints, i.e., unconstrained, semi-constrained with the root condition, and fully-constrained with both root and consistency conditions. We focus on the learning of 3-step linear multistep methods, which we subsequently applied to solve three model PDEs, i.e., the one-dimensional heat equation, the one-dimensional wave equation, and the one-dimensional Burgers' equation. The results show that the prediction error of the learned fully-constrained scheme is close to that of the Runge-Kutta method and Adams-Bashforth method. Compared to the traditional methods, the learned unconstrained and semi-constrained schemes significantly reduce the prediction error on coarse grids. On a grid that is 4$\times$ coarser than the reference grid, the mean square error shows a reduction of up to an order of magnitude for some of the heat equation cases, and a substantial improvement in phase prediction for the wave equation. On a 32$\times$ coarser grid, the mean square error for the Burgers' equation can be reduced by up to 35\% to 40\%.

\keywords
{Time integration scheme, Partial differential equation, Numerical simulation on coarse grids, Neural networks, Mathematical constraints}

\section{Introduction}\label{sec:Introduction}

Many environmental and engineering problems, encompassing a broad spectrum of spatial and temporal scales, are described by partial differential equations (PDEs) of high dimension. Since resolving all the scales often demands exceedingly extensive computational resources, solving the high-dimensional PDEs poses a great challenge to the current computing systems. For example, direct numerical simulations of high-Reynolds number turbulent flows can be particularly challenging in this regard.
Solving the PDEs on coarse grids, reducing the dimension of the problem, can substantially reduce the computational costs, while meanwhile introduces discretization errors. Discretization errors depend not only on the spatial discretization schemes but also on the time integration schemes. In this study, we propose to learn time integration schemes to reduce the {\color{black}error of solving PDEs on coarse grids}. 

The traditional approach to reduce the {\color{black}prediction error} on coarse grid is to directly account for the effect of the physics of the unresolved scales. 
One way is to solve the spatially filtered PDEs instead of the original PDE. However, the filtering procedure on the nonlinear term often introduces a new unclosed term, i.e., the subgrid term. 
The subgrid term can be modelled explicitly by establishing its relation with the resolved scales. For instance, the eddy viscosity model for large-eddy simulation (LES) of turbulent flows governed by the Navier-Stokes (NS) equations~\cite{smagorinsky1963general, germano1991dynamic}. The subgrid term can also be modelled implicitly by changing the properties of the schemes for discretizing the spatial derivatives. One example is the Implicit Large-Eddy Simulation (ILES) method for solving turbulent flows~\cite{grinstein2007implicit}. The success of such approaches largely depends on our understanding of the physics of the unresolved scales. Moreover, the developed models are often valid for only a range of unresolved scales. When the dominant physics of the unresolved scales change from one to the other as the cut-off scale changes, a grey region appears where neither of the models in each regime works. For instance, the grey region in the microscale and mesoscale simulations of the atmospheric flow in the meteorology research~\cite{Wyngaard2004, Honnert2020}.  

The machine learning methods are revolutionizing the numerical methods for solving PDEs~\cite{Raissi2017a, Sirignano2018, Lu2021, Wang2021, Zhang2019, ren2022phycrnet}. 
Such approaches approximate the solution based on machine learning models~\cite{Karniadakis2021, Cuomo2022}, e.g., neural networks~\cite{Poggio2017}, avoiding the need for a grid to discretize the computational domain. In the review paper by \cite{Blechschmidt2021}, {\color{black}three types of the methods based on neural networks were reviewed}: 1) the physics-informed neural networks (PINN)~\cite{Raissi2017a, Raissi2017b, Raissi2019}, in which the PDEs with initial and boundary conditions are approximately enforced via loss functions; 2) the methods based on the Feynman–Kac formula~\cite{Beck2018}, which approximate the solution of a PDE as the expectation of a stochastic process; and 3) the methods based on the solution of backward stochastic differential equations~\cite{E2017, Han2018}, in which deep neural networks (DNN) are employed for computing the gradient of solutions. As noted in the review~\cite{Blechschmidt2021}, the last two methods showed promising results for high-dimensional linear and semilinear systems. The methods based on PINN are capable of handling complex nonlinear PDEs and inverse problems with incomplete models and imperfect data, while they are not competitive for solving well-posed, high-dimensional forward problems when compared with the traditional grid-based numerical methods~\cite{Karniadakis2021, Blechschmidt2021}. 

Improving the predictive capability of grid-based numerical methods using the machine learning methods received lots of attention as well~\cite{duraisamy2019turbulence, zhu2019machine,zhou2019subgrid, meng2022artificial, Yang_etal_PRF_2019, zhou2021wall, Zhang2022, Zhang2023, bezgin2021data, kossaczka2022neural, 2018Learning}.
Some efforts are particularly focused on improving the traditional subgrid models for solving PDEs on coarse grids, such as models for turbulent flow simulations~\cite{duraisamy2019turbulence, zhu2019machine,zhou2019subgrid, meng2022artificial, Yang_etal_PRF_2019, zhou2021wall, Zhang2022, Zhang2023}. 
Other attempts were made to improve the performance of the spatial discretization schemes for PDEs~\cite{bezgin2021data, kossaczka2022neural, 2018Learning}. 
Particularly, the spatial discretization schemes for solving PDEs on coarse grids were learned using neural networks in the work by Bar-Sinai et al.~\cite{2018Learning}. The results of the one-dimensional Burgers' equation showed that the learned discretization schemes can reduce the prediction {\color{black}error} and extend the median valid simulation time for very coarse grids.    

The stability and accuracy of solving PDEs on coarse grids depend on both temporal and spatial discretization schemes. Empowering time integration schemes through machine learning methods, however, has been overlooked in existing studies. When solving PDEs on coarse grids, errors in the solution due to coarse spatial resolution accumulate over time. For example, this can lead to phase errors when solving the wave equation on coarse grids. In this study, we employ neural networks to learn time integration schemes to reduce the {\color{black}error} of solving PDEs on coarse grids. Specifically, we propose three learning approaches with different constraints and test the learned schemes using the one-dimensional heat equation, wave equation, and Burgers' equation.

In the rest of the paper, the three learning approaches are first described and applied to learn a 3-step linear multistep method in section~\ref{sec:model}. Then, the test results for the one-dimension heat equation, wave equation, and Burgers' equation are presented sequentially in section~\ref{sec:results}. Lastly, the conclusions are drawn in section~\ref{sec:conclusions}. 

\section{Learning of time integration schemes}
\label{sec:model}
In this section, we describe the approach for developing data-driven time integration schemes. The 3-step linear multistep method is taken as an example. Similar ideas can be applied to other time integration schemes. The general formulation of the linear multistep method is given in section~\ref{sec:linear multistep} focusing on the constraints for the stability and consistency of the method. The three different models for imposing the constraints are then presented in section \ref{sec:constraints}. The definition of the loss function is given in section~\ref{sec:loss}. Lastly, the details for the application in a 3-step linear multistep method is given in section~\ref{sec:third-order method}.
\subsection{General formulation of the linear multistep method} \label{sec:linear multistep}
We consider the following ordinary differential equation (ODE) with independent variable $t$
\begin{linenomath*}
\begin{equation}
\frac{\partial v}{\partial t} = F(x,t),
\end{equation}
\end{linenomath*}
where $x$ is a parameter (which can be the spatial coordinate). The linear multistep method of $k$-step for discretizing the above ODE  
can be expressed as
\begin{linenomath*}
\begin{equation}
    \label{eq: linear multistep}
\alpha_{k} v^{n+1}= -\sum_{i=0}^{k-1} {\alpha_{i} v^{n+i-k+1}} +\Delta t \sum_{i=0}^{k} {\beta_{i} F_{n+i-k+1}},
\end{equation}
\end{linenomath*}
where $\alpha_{i}$ and $\beta_{i}$ are real coefficients, $v^{n}$ represent the approximate value of $v(t_{n})$ at $t=t_{n}$, $\Delta t$ is the size of time step and $F_{n+i} = F(x,t_{n+i})$, $t_{n+i}=t_{n}+i\Delta t$. For the sake of simplicity, we assume that $\alpha_{k}=1$, $\alpha_{0}^{2}+\beta_{0}^{2} \neq 0$. 

The generating polynomials of the multistep method, which can be written as
\begin{linenomath*}
\begin{align}
\rho (\chi)&=\alpha _{k}\chi ^{k}+\alpha _{k-1}\chi ^{k-1}+\cdots +\alpha _{0}, \label{first_poly}\\
\sigma (\chi)&=\beta _{k}\chi ^{k}+\beta _{k-1}\chi ^{k-1}+\cdots +\beta _{0}, \label{second_poly}
\end{align}
\end{linenomath*}
are employed to obtain the consistency and stability constraints~\cite{hairer1993solving}. The obtained consistency condition is 
\begin{linenomath*}
\begin{equation}
    \label{eq:consistency}
    \rho (1)=0,\quad\rho' (1)=\sigma (1),   
\end{equation}
\end{linenomath*}
which is also the condition that the linear multistep method should satisfy to be of first-order accuracy \footnote{According to Taylor expansion, the necessary and sufficient condition for a linear multistep method of $k$-step to be of order $P$ is $\sum_{i=0}^{k} \alpha_{i} =0$, $\sum_{i=0}^{k} \alpha_{i}i^{s}=s\sum_{i=0}^{k} \beta_{i}i^{s-1}, s=1,2,\cdots P.$}.

The stability condition we employ here is called zero-stability, that a general method is stable when $\Delta t\to 0$. For a linear multistep method, it is stable if $\rho (\chi )$ satisfies the root condition~\cite{hairer1993solving}, i.e., the moduli of the roots of $\rho (\chi )$ are less than or equal to 1, and if one root’s modulus is 1, the root must be single. That is to say, the roots of $\rho (\chi )$ 
are on or within the unit circle, and the roots on the unit circle are single. It is worth mentioning that the convergence and stability of the linear multistep method are equivalent under the premise of the consistency condition.

\subsection{Three models for learning time integration schemes}
\label{sec:constraints}

The stability and consistency constraints discussed in the last subsection are imposed during the training of the neural network model, with the aim to increase the interpretability and stability of the learned time integration schemes. Three different strategies for applying the constraints are employed, i.e., 1) the BP (backpropagation) neural network model with no constraints (the unconstrained model, Un-con); 2) the BP neural network model with root condition (stability condition) enforced (the semi-constrained model, Semi-con); 3) the BP neural network model with both root condition and consistency condition enforced (the fully-constrained model, Full-con). In the following, the root condition and consistency condition will be reformulated in a way to facilitate their implementation as constraints for training the neural network models. As the results from the explicit Runge-Kutta method of order 3(2) with adaptive time-stepping will be employed as the reference for examining the accuracy, in the following we focus on deriving the constraints for the data-driven explicit linear multistep method of 3-step ($k=3$ in Eq.~(\ref{eq: linear multistep})).

The consistency condition, which is equivalent to imposing the first-order constraint to a third-order explicit linear multistep method, is implemented in a way that some parameters are the direct outputs of the neural network with the rest of the parameters (i.e., $\alpha_{i}$ and $\beta_{i}$ in Eq.~(\ref{eq: linear multistep})) adjusted to satisfy the constraint (i.e., Eq.~(\ref{eq:consistency})). 

In regard to the root condition (stability condition), i.e., the moduli of the roots of $\rho (\chi )$ are less than or equal to 1 ($\left \|\chi  \right \|  \le 1$), and if one root’s modulus is 1, the root must be single, we employ the transform $\chi = \frac{1+z}{1-z}$ to map the above constraint to the one that the real parts of all roots of a real coefficient polynomial with respect to $z$ are negative. With $Re(z)<0$, the original root condition is replaced by a stricter one, $\left \|\chi  \right \|  < 1$. Therefore, the $\rho (\chi )$ with the root condition promised can be turned into a Hurwitz Polynomial $ \psi (z)$~\cite{gantmacher2005applications}. Specifically, the Hurwitz Polynomials are in the following form, 
\begin{linenomath*}
\begin{equation}
    \psi_2 (z)=(1-z)^2\rho_{2}\left(\frac{1+z}{1-z}\right)=(1-\alpha_{1}+\alpha_{0})z^2+2(1-\alpha_{0})z+(1+\alpha_{1}+\alpha_{0}),
\end{equation}
\end{linenomath*}
\begin{linenomath*}
\begin{equation}
    \psi_3 (z)=(1-z)^3 \rho_{3}\left(\frac{1+z}{1-z}\right)=(1-\alpha_{2}+\alpha_{1}-\alpha_{0})z^3+(3-\alpha_{2}-\alpha_{1}+3\alpha_{0})z^2+(3+\alpha_{2}-\alpha_{1}-3\alpha_{0})z+(1+\alpha_{2}+\alpha_{1}+\alpha_{0}),
\end{equation}
\end{linenomath*}
which correspond to the following quadratic and cubic $\rho (\chi )$ polynomials $\rho_{2} (\chi )=\chi ^2+\alpha_{1}\chi+\alpha_{0}$ and $\rho_{3} (\chi )=\chi ^3+\alpha_{2}\chi ^2+\alpha_{1}\chi+\alpha_{0}$, respectively.  

According to the Routh-Hurwitz criterion~\cite{gantmacher2005applications}, the constraint that the real parts of the roots of $\psi_2 (z)$ and $\psi_3 (z)$ are negative is equivalent to the condition that their polynomial coefficients satisfy the following inequalities,
\begin{linenomath*}
\begin{equation}
\label{eq:coeff}
\begin{cases}
1-\alpha_{1}+\alpha_{0}>0, \\
1-\alpha_{0}>0, \qquad \quad \text{(a)}\\
1+\alpha_{1}+\alpha_{0}>0,\\
\end{cases}\qquad
\begin{cases}
1-\alpha_{2}+\alpha_{1}-\alpha_{0}>0, \\
1+\alpha_{2}+\alpha_{1}+\alpha_{0}>0,\\
1-\alpha_{1}+\alpha_{2}\alpha_{0}-\alpha_{0}^2>0, \qquad \text{(b)}\\
1-\alpha_{0}>0,\\
1+\alpha_{0}>0.\\
\end{cases}
\end{equation}
\end{linenomath*}
Therefore, the reformulated root condition (i.e., Eq.~(\ref{eq:coeff})) is obtained. 

For the unconstrained model, no constraints are imposed. For the semi-constrained model, the root condition constraint, i.e., Eq.~(\ref{eq:coeff}b), is enforced. In term of the fully-constrained model, the obtained coefficients of the linear multistep method are required to meet the consistency condition and the root condition. As the consistency condition (Eq.~(\ref{eq:consistency})) implies that $\rho(\chi)$ has a root of 1, the $\rho_3(\chi)$ polynomial can be written as
\begin{linenomath*}
\begin{equation}
\label{eq:pq 2-ord}
\rho_3(\chi)=(\chi-1)(\chi^2+p\chi+q).
\end{equation}
\end{linenomath*}
The root condition is then reduced to that
the rest of the two roots of $\rho_2(\chi)$ must lie within the unit circle. The coefficients of $\chi^2+p\chi+q$ should meet the Routh-Hurwitz criterion, i.e., Eq.~(\ref{eq:coeff}a). With $p$, $q$ given by the learned neural network model, the values of $\alpha_{0}$, $\alpha_{1}$ and $\alpha_{2}$ can then be obtained. 

\subsection{Loss function}
\label{sec:loss}
The loss function consists of two parts, i.e., the part related to the {\color{black}error of the solution predicted by the learned time integration scheme}
, and the other part for the constraints imposed on the scheme coefficients given by the neural network model. The first part of the loss function ($\text{MSE}(v_{pred},v_{true})$) is defined as the mean squared error between the predictions on the coarse grid ($v_{pred}$) and those on the same grid produced by coarsening the solutions from the high resolution grid ($v_{true}$). 
The second part of the loss function is in the form of the internal penalty function, which originates from the inequalities of the root condition constraint. The internal penalty function, which is also called the barrier function, establishes a barrier at the boundary of the feasible region to prevent the iteration from leaving the area~\cite{luenberger1984linear}. When the output values of the neural network are close to the boundary of the feasible region, the value of the internal penalty function tends to become infinity. In this work, the barrier function of the following reciprocal form, 
\begin{linenomath*}
\begin{equation}
B_{1}(\vec{\alpha})=\sum_{i} \frac{1}{\left |   g_{i}(\vec{\alpha})\right |+\varepsilon },\quad \text{(a)} \quad \quad  \quad \quad B_{2}(p, q)=\sum_{i} \frac{1}{g_{i}(p, q)},\quad \text{(b)}
\label{eq:barrier}
\end{equation}
\end{linenomath*}
that the boundary of the feasible region is determined by $ g_{i}(\cdot)>0$, is employed. The barrier function for the semi-constrained model Eq.~(\ref{eq:barrier}a) sets the constraints on the values of $\alpha_{i}$ directly. Because the initial output value may be at the boundary of the feasible region, the application of absolute values and $\varepsilon$ here is to prevent the denominator from changing to $0$ during training, resulting in subsequent training not being possible. In addition we need to test the output of the semi-constrained model after each training. If the output coefficients are not in the feasible domain, the model will be retrained with modified initial weights.

For the fully-constrained model, the barrier function Eq.~(\ref{eq:barrier}b) is enforced on $p$, $q$ as shown in Eq.~(\ref{eq:pq 2-ord}), in which the consistency condition $\rho (1)=0$ is guaranteed.  Due to the initial output of the fully constrained model being within the feasible domain, it will not run out of the feasible domain at small learning rates. There is no need for techniques to ensure that the denominator is not zero. The specific initial network settings are shown in subsection~\ref{sec:third-order method} and section~\ref{sec:results}.
As a consequence, the loss function can be written as
\begin{linenomath*}
\begin{equation}
\text{loss}=\text{MSE}(v_{pred},v_{true})+\gamma B,
\label{eq:loss}
\end{equation}
\end{linenomath*}
where $\gamma$ is an adjustable hyperparameter and set to 0 for the unconstrained model. The $\varepsilon$ in Eq.~(\ref{eq:barrier}a) is a small quantity and should be adjusted with $\gamma$ changes. 

\subsection{Learning of 3-step linear multistep methods}\label{sec:third-order method}
The data preparation and model training of a 3-step linear multistep method are presented to show the procedure of developing data-driven time integration schemes based on the BP neural networks, which is then applied to the one-dimensional heat equation, wave equation and Burgers' equation.

The outputs of the neural network model are the coefficients $\alpha_{i}$ and $\beta_{i}$ for the linear multistep formulation (Eq.~(\ref{eq: linear multistep})). The inputs consist of coarse-grained solutions from the previous time step, i.e., $v^{n}$ at different spatial locations (with the solutions $v^{n-1}$, $v^{n-2}$ for computing the right-hand-side term for the 3-step linear multistep method). Consequently, the learned coefficients are influenced by spatial and temporal variations in the solutions. Since the internal penalty function is employed, the initial values of the outputs of the neural network need to be specified in the feasible region.
A small learning rate is employed, such that the step size of updating the hyperparameter of network is small, which ensures the outputs are located in the feasible region. Max-min normalization is used for inputs during training and testing to ensure valid constraints during testing.

To train the data-driven model for a 3-step linear multistep method, the coefficients of the third-order explicit Adams methods, which can be written as
\begin{linenomath*}
\begin{equation}
\label{3-ordadams}
v^{n+1}=v^{n}+\Delta t\left(\frac{23}{12} F_{n}-\frac{16}{12} F_{n-1}+\frac{5}{12} F_{n-2}\right), 
\end{equation}
\end{linenomath*}
are employed to set the initial output of the neural network. To prevent the training results from falling into local optimization, small perturbations are added to the coefficients. With the assumption that the optimal coefficients are close to the initial guess, the initial biases of the network are set manually, and the initial weights are set as small random numbers for realizing the optimization in small step sizes and ensuring that the output coefficients are still in the feasible region when the inputs are not in the training set. 

\begin{figure}[htb!]
\centering
\includegraphics[width=14cm]{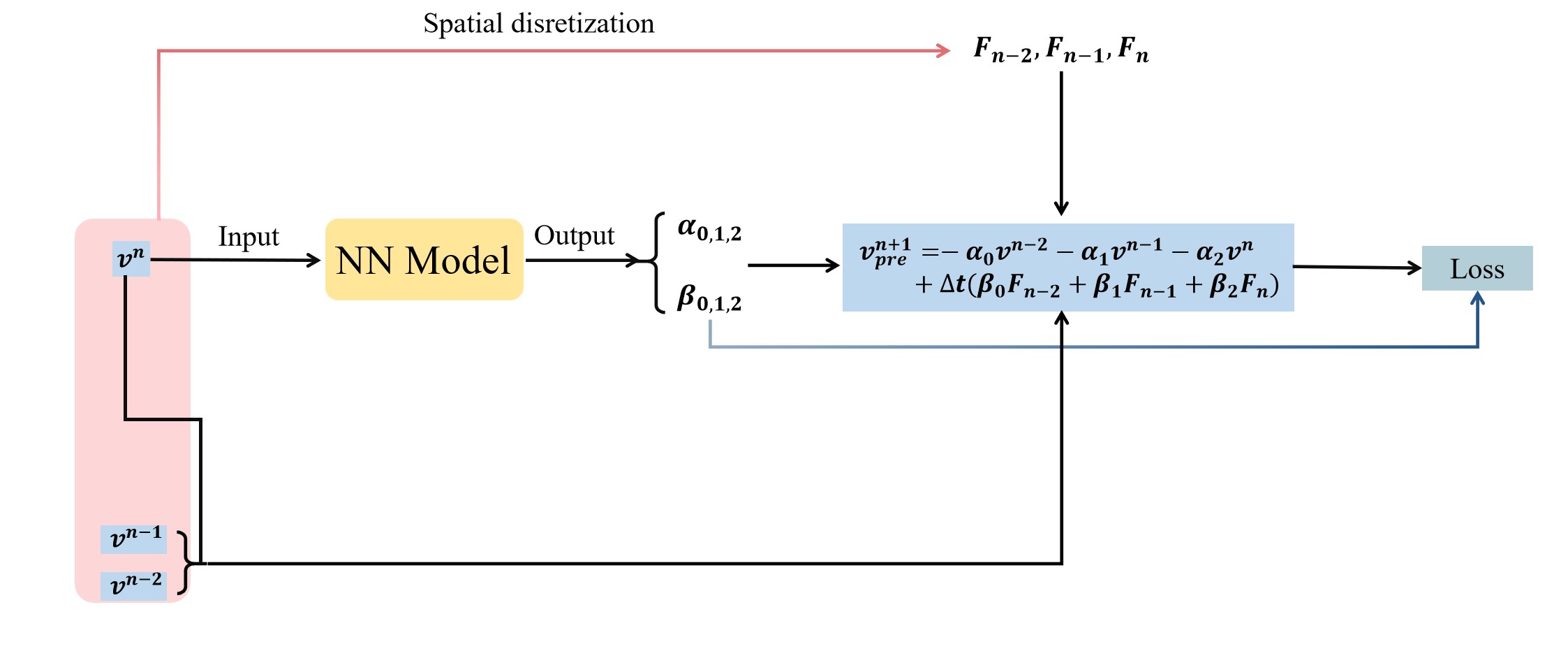}
	  \caption{ A schematic for the training procedure of a data-driven 3-step linear multistep scheme. }\label{fig:training}
\end{figure}

A schematic showing the training process is shown in Figure \ref{fig:training}. During each iteration of the model training, the solution from the previous time step is first given as the input. With the input, the scheme coefficients are then given by the NN model (neural network model). At last, the loss is computed as the mean square error of the predicted solution $v_{pre}^{n+1}$ and the constraints for the coefficients (i.e., Eq.~(\ref{eq:loss})), in which $v_{pre}^{n+1}$ is computed by advancing the equation for one step using the output coefficients of the NN model, the solution at previous time steps, and the right-hand-side term. 

The employed neural network is composed of three layers, i.e., the input layer with the number of neurons equals to the number of grid points in $x$ , a hidden layer with 20 neurons, and the output layer with 6 neurons for the unconstrained model and the semi-constrained model, and 4 neurons for the fully-constrained model. Due to the initial loss being below $10^{-5}$, small learning rates, somewhere around $10^{-7}$, are employed to adjust the initial guess of the coefficients. The three models are all trained using the Adam optimizer. The activation function is ReLU.
As the number of biases in the first two layers is greater than the output, the values of the additional bias will default to the last number of the manually set bias for the output layer. The models are trained in the framework of TensorFlow 1.0.

\begin{figure}[htb!]
\centering
\includegraphics[width=14cm]{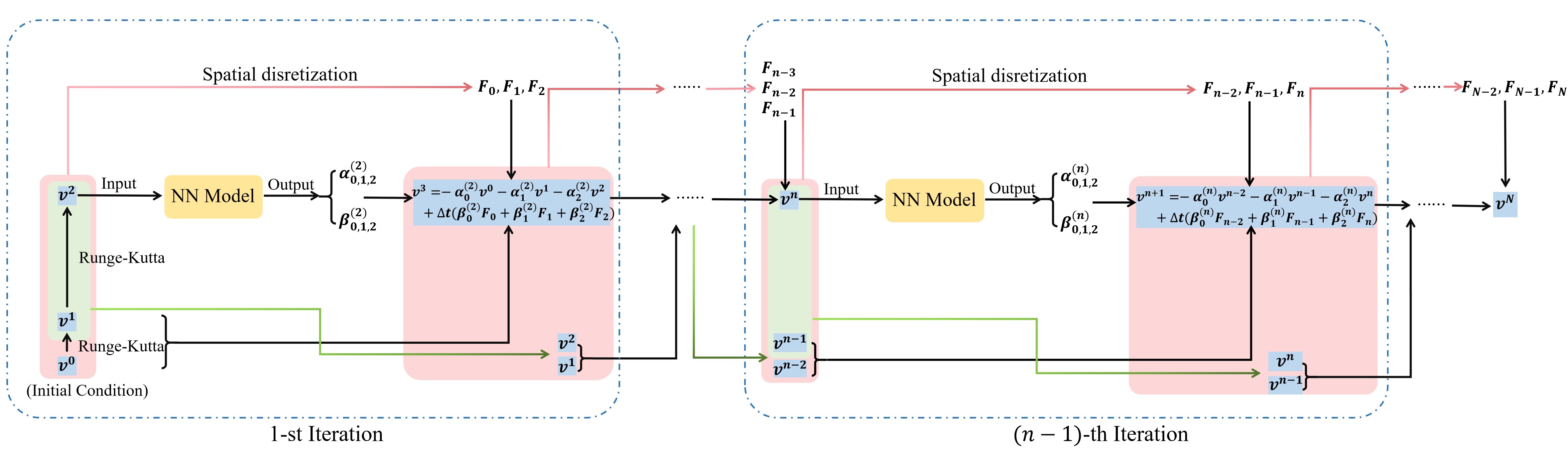}
	  \caption{A schematic for the procedure of using the learned 3-step linear multistep scheme.}\label{fig:calculate}
\end{figure}

The training data are generated by coarsening the numerical solution on the fine grid using the cell-averaging approach. The spatial discretizations based on the finite volume method are performed, which has good stability and conservation properties, being consistent with the cell-averaging approach employed for coarsening. 

The procedure for using the learned data-driven 3-step time integration scheme is shown in Figure~\ref{fig:calculate}. In this figure, the process of advancing in time is shown by the green arrows. At beginning, the solutions $v^{1}$ and $v^{2}$ are calculated from the initial solution $v^{0}$ using the Runge-Kutta method. The black arrows indicate the process of employing the data-driven time integration scheme, while the red arrows indicate the calculation of the right-hand-side terms.

\section{Results} \label{sec:results}
In this section, we present the results from the data-driven time integration schemes and compare them with those from the Runge-Kutta method of order 3(2) with adaptive time-stepping. In subsection~\ref{sec:Heat Equation}, we apply the data-driven time integration schemes with the finite volume method to solve 1-D heat equations with different thermal diffusivities on a grid $4\times$ coarser than the reference grid. And then we consider first-order wave equations with varying wave speeds in subsection~\ref{sec:Wave Equation}, and analyze why data-driven time integration schemes can achieve highly accurate results. 
In subsection~\ref{sec:burgers}, we test the data-driven temporal schemes on the Burgers' equation on a grid 32$\times$ coarser than the reference grid, with the spatial derivatives approximated using the data-driven schemes proposed in~\cite{2018Learning}. 
 
\subsection{1-D Heat Equation} 
\label{sec:Heat Equation}
The heat equation considered here is 
\begin{linenomath*}
\begin{equation}
\begin{cases}
\frac{\partial v}{\partial t}=  \frac{\partial}{\partial x} \left(\lambda\frac{\partial v}{\partial x}\right), \quad 0 \leqslant x \leqslant 1, \; 0 \leqslant t \leqslant 1,\\
v(x, 0)=\sin(2 \pi x), \quad 0 \leqslant x \leqslant 1, \\
\end{cases} \label{eq:heat_eq}
\end{equation}
\end{linenomath*}
where $\lambda$ is thermal diffusivity which expresses the ability of the a system tending to a uniform temperature in heating or cooling~\cite{anderson2016computational}. This equation describes the heat conduction or diffusion in one-dimensional isotropic medium. Eq.~(\ref{eq:heat_eq}) employs the periodic boundary condition with the domain of size $L=1$. The exact solution is
\begin{linenomath*}
\begin{equation}
v(x, t)=\mathrm{e}^{-4\pi^{2}\lambda t} \sin(2 \pi x ).\label{eq:heat_exact}
\end{equation}
\end{linenomath*}

The thermal diffusivity $\lambda$ here is set in the range of 0.1 to 1.  We apply the data-driven time integration schemes trained under a particular thermal diffusivity to other thermal diffusivities in this range. The first-order finite-volume method is employed for spatial discretization during both training and testing.

The three different data-driven time integration schemes with/without constraints (i.e., the unconstrained, the semi-constrained, and the fully-constrained schemes) are trained on a training set with a fine grid solution coarsen by a factor of 4 for $\lambda= 0.5$, and then tested on  $\lambda \in \left \{ 0.1, 0.2, 0.3, 0.4, 0.6, 0.7, 0.8, 0.9, 1.0 \right \} $. We assume that, for example, if time integration schemes have lower errors at $\lambda=0.7$ and 0.8, then it is highly likely that schemes will perform better in $\lambda \in \left[0.7,0.8\right] $. The fine grid solution has 64 grid points in the $x$ direction for $0 \leqslant x \leqslant 1$, which means there are 16 grid points on the 4$\times$ coarse grid. The size of the time step is set to 0.0001.

Specifically, to learn the data-driven 3-step linear multistep model, the initial weights are uniformly distributed in the range from $-0.0005$ to $0.0005$ for the unconstrained model and the semi-constrained model, with the same initial biases $[0, 0, -1, 5/12, -4/3, 23/12]$ for $\alpha_{0,1,2}$ and $\beta_{0,1,2}$. With the consistency condition enforced, the number of outputs in the fully-constrained model is 4, i.e., $p$, $q$, $\beta_{0}$ and $\beta_{1}$. Considering that $p, q=0$ for the generating polynomial $\rho(\lambda)$ of Eq.~(\ref{3-ordadams}) when simplifying it as the form of Eq.~(\ref{eq:pq 2-ord}), the initial biases of the fully-constrained model are set as $[0, 0, 5/12, -4/3]$, and the initial weights are uniformly distributed in the range from $-0.005$ to $0.005$. The learning rate is $10^{-7}$ for both unconstrained model and semi-constrained model and is $5\times10^{-7}$ for fully-constrained model. The three models are trained by Adam optimizer for 8000 steps. The $\gamma$ in loss function Eq.~(\ref{eq:loss}) is set to $10^{-18}$ for the semi-constrained model and $10^{-12}$ for the fully-constrained model. {\color{black}The purpose is that when the output coefficients are in the feasible domain, the penalty function will minimize its impact on the main part of the loss, i.e. the prediction error.}

\begin{figure}[htb!]
\centering
\includegraphics[width=14cm]{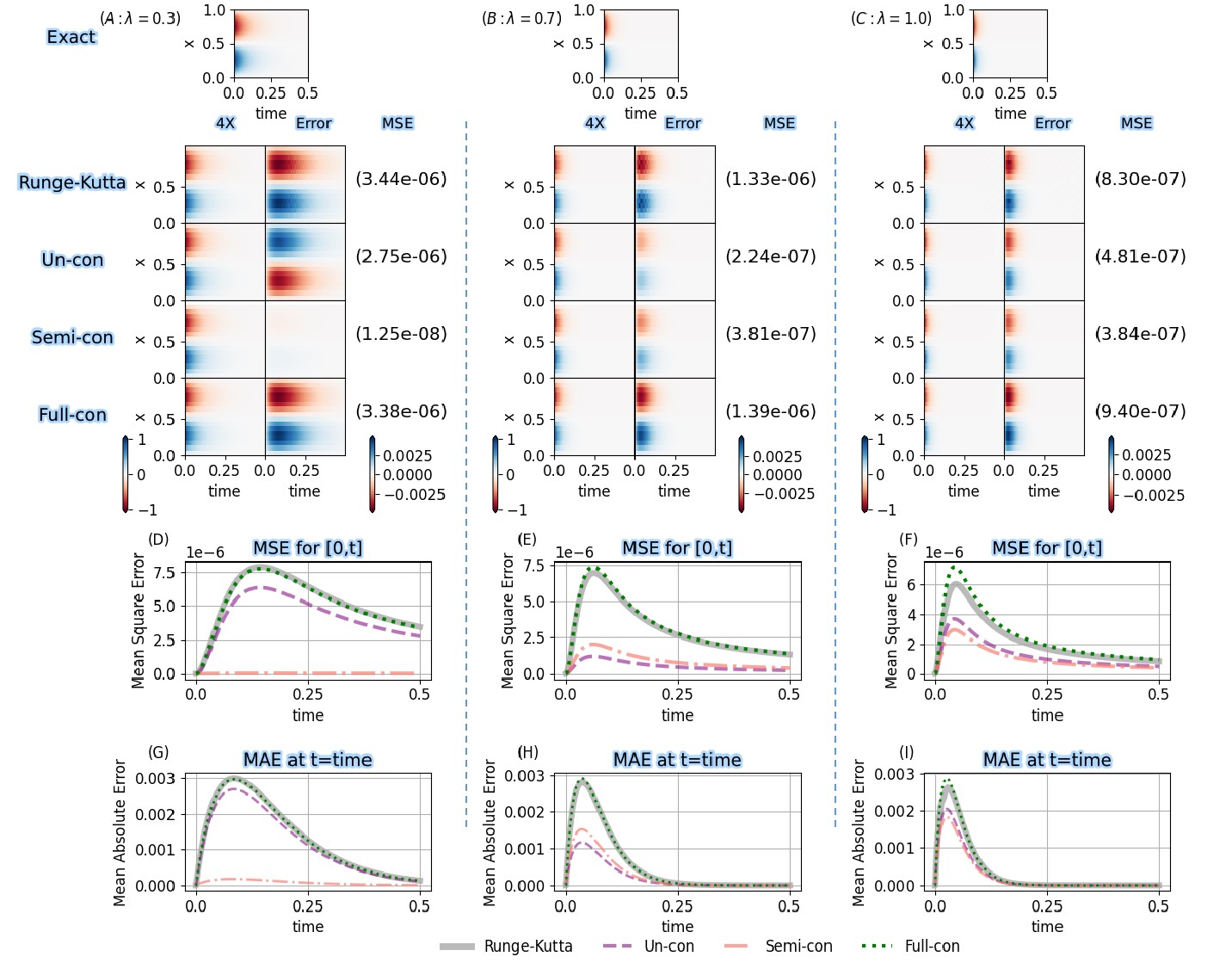}
\caption{Test results of the 1-D heat equation for (A) $\lambda=0.3$, (B) $\lambda=0.7$ and (C) $\lambda=1.0$ and their curves of error over time (D, E, F, G, H, I). The subgraph on the top of (A, B, C) is the exact solution obtained from Eq.~(\ref{eq:heat_exact}) and then followed by realizations of solutions and absolute error distribution from different time integration schemes. The numbers in brackets on the right of (A, B, C) are the mean square error averaged over the whole domain. The error shown in (D, E, F) is the mean squared error obtained by averaging the error in space ($[0, 1]$) and time ($[0,t]$). The error shown in (G, H, I) is the mean absolute error averaged over space ($[0, 1]$) at $t$ instant.
(D, G) is the error curves for (A): $\lambda=0.3$ and (E, H) for (B): $\lambda=0.7$, (F, I) for (C): $\lambda=1.0$.
}
\label{fig:heat_equation}
\end{figure}

Figure \ref{fig:heat_equation} shows the results of the 1-D heat equation for $\lambda=0.3, 0.7$ and 1.0. Since the solution tends to be a constant after a time greater than 0.5, only the solution and error for t from 0 to 0.5 are plotted here. The results for a larger range of t and other $\lambda$ can be found in table \ref{tab:heat_error} in the Appendix. From Figure \ref{fig:heat_equation} and table \ref{tab:heat_error}, it can be seen that unconstrained model and semi-constrained model have lower values of mean square error (MSE) and mean absolute error (MAE) than the Runge-Kutta method. Little difference is observed between prediction from the fully-constrained model and Runge-Kutta method. We believe this is because the fully-constrained model needs to find coefficients satisfying both the root condition and the consistency condition, leaving almost no space for optimisation.

As can be seen from the table \ref{tab:heat_error} and figure \ref{fig:heat_equation}, while the unconstrained model  {\color{black} works well} for $\lambda=0.3$ to 1 (especially for $\lambda$ close to the training set around $0.5$), it performs the worst for $\lambda=0.1$ and 0.2, with the MSE one to two orders of magnitude greater than the Runge-Kutta method. This is understandable as the unconstrained model has the space to learn the mechanism for error reduction, but with a narrow scope of generalisation as a result of not satisfying certain mathematical constraints. 

In contrast, the generalisation ability of the semi-constrained model is well demonstrated, with its predictions closer to the exact solution when compared with the Runge-Kutta method for the cases with $\lambda$ from 0.2 to 1.  
This can be explained by the fact that the root constraint enhances the stability of the learned scheme, and meanwhile leaves room for optimisation. 

\subsection{1-D Wave Equation} \label{sec:Wave Equation}
In this subsection one-dimensional first-order wave equation (i.e. linear convection equation~\cite{anderson2016computational}) with wave speed $c$ is considered, which can be written as
\begin{linenomath*}
\begin{equation}
\begin{cases}
\frac{\partial v}{\partial t}+ c\frac{\partial v}{\partial x}=0, \quad 0 \leqslant x \leqslant 1, \; 0 \leqslant t \leqslant 1,\\
v(x, 0)=\sin(4 \pi x), \quad 0 \leqslant x \leqslant 1.\\
\end{cases} \label{eq:wave_eq}
\end{equation}
\end{linenomath*}
where the wave speed $c > 0$. The wave equation describes the propagation of a wave in a homogeneous medium with a velocity $c$ in the $x$ direction. Eq.~(\ref{eq:wave_eq}) employs the periodic boundary condition with the domain of size $L=1$. The exact solution is
\begin{linenomath*}
\begin{equation}
v(x, t)=\sin \left [ 4 \pi (x - ct) \right ] .\label{eq:wave_exact}
\end{equation}
\end{linenomath*}

The wave speed $c$ here is set from 0.1 to 1. The spatial discretization scheme is the {\color{black}second order} finite-volume method. The training set is generated by coarsening the numerical solution on the fine grid with $c=0.5$. The test set has the values of $c\in \left \{0.1, 0.2, 0.3, 0.4, 0.6, 0.7, 0.8, 0.9, 1.0 \right \}$. The fine grid case has points of 64 in $x$-direction. The cell-averaging method is employed for coarsening find grid solution to obtain the training data. The size of the time step is set to 0.0001.

To train the model, the initial weights are set uniformly distributed in the range of $\left [-0.0005, 0.0005 \right ]$ for the unconstrained model and $\left[0,0.0005\right ]$ for the semi-constrained model while $\left [-0.005, 0.005 \right ]$ for the fully-constrained model. The initial biases, learning rate, and optimizer are identical with the 1-D heat equation case but training for 10000 steps.

\begin{figure}[htb!]
\centering
\includegraphics[width=14cm]{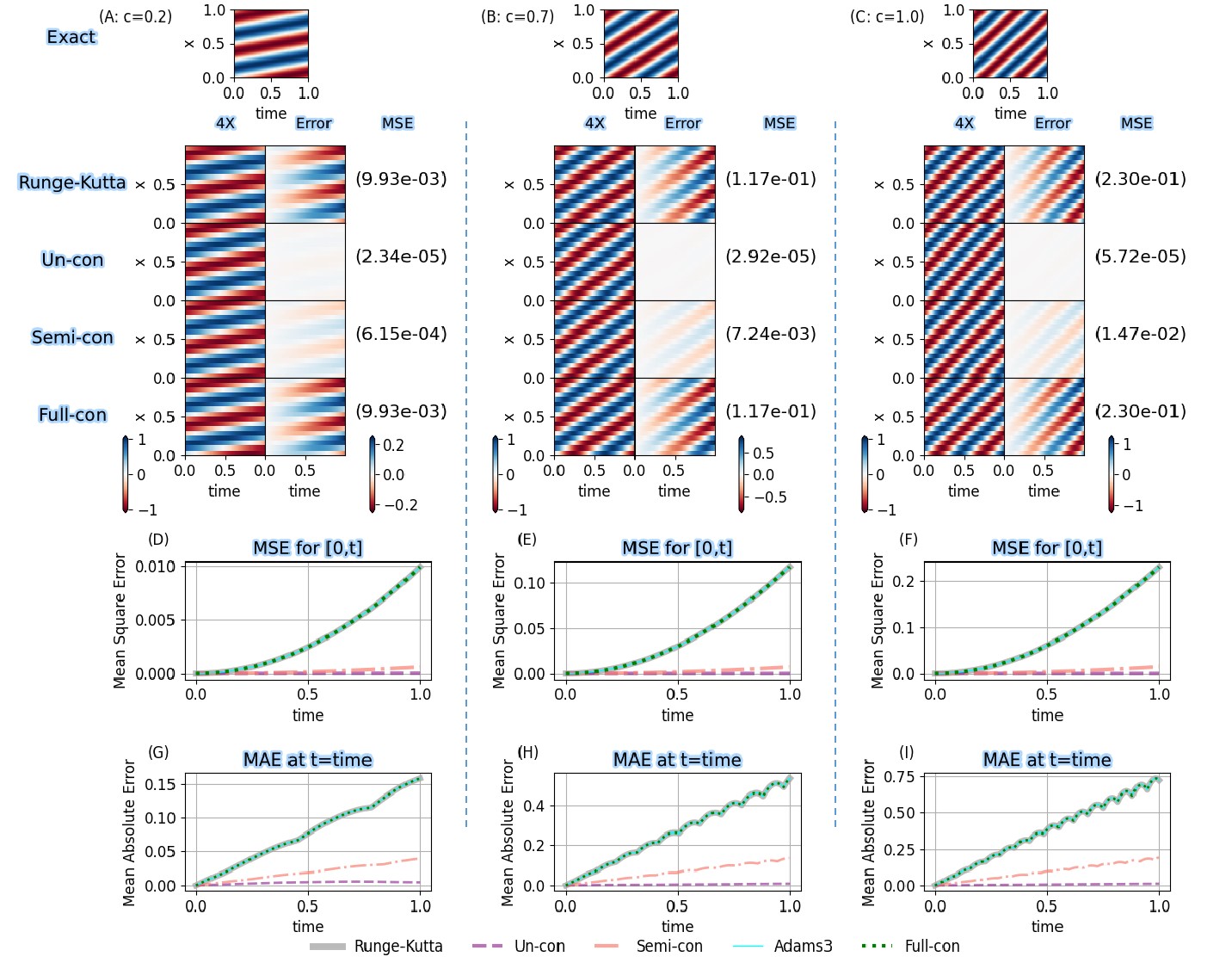}
\caption{Test results of wave equation for (A) $c=0.2$, (B) $c=0.7$ and (C) $c=1.0$ and their curves of error over time (D, E, F, G, H, I). The subgraph on the top of (A, B, C) is the exact solution obtained from Eq.~(\ref{eq:wave_exact}) and then followed by realizations of solutions and absolute error distribution from different time integration schemes. The numbers in brackets on the right of (A, B, C) are the mean square error averaged over the whole domain. The error shown in (D, E, F) is the mean squared error obtained by averaging the error in space ($[0, 1]$) and time ($[0,t]$). The error shown in (G, H, I) is the mean absolute error averaged over space ($[0, 1]$) at $t$ instant. (D, G) is the error curves for (A): $c=0.2$ and (E, H) for (B): $c=0.7$, (F, I) for (C): $c=1.0$.
}
\label{fig:wave_equation}
\end{figure}

Figure \ref{fig:wave_equation} shows the results of the 1-D wave equation for $c=0.2, 0.7$ and 1.0. The MSE and the MAE for other $c$ could be found in table \ref{tab:wave_error} in the Appendix. Here we clearly see that the coarse-grained prediction of the unconstrained model for $c$ from 0.1 to 1 is the closest to the exact solution among the considered methods, with the values of MSE and MAE reduced by at least one order of magnitude in comparison with the Runge-Kutta method and the third-order Adams methods. Significant and consistent improvements are also observed for the semi-constrained model, with at least about 90 percent reduction in MSE and 60 percent reduction in MAE for each $c$. The fully-constrained model, on the other hand, predicts almost the same solutions as the Runge-Kutta method and the third-order Adams method. For the sake of comparison, the third-order Adams method will be used here as the baseline method.

\begin{figure}[htb!]
\centering
\includegraphics[width=14cm]{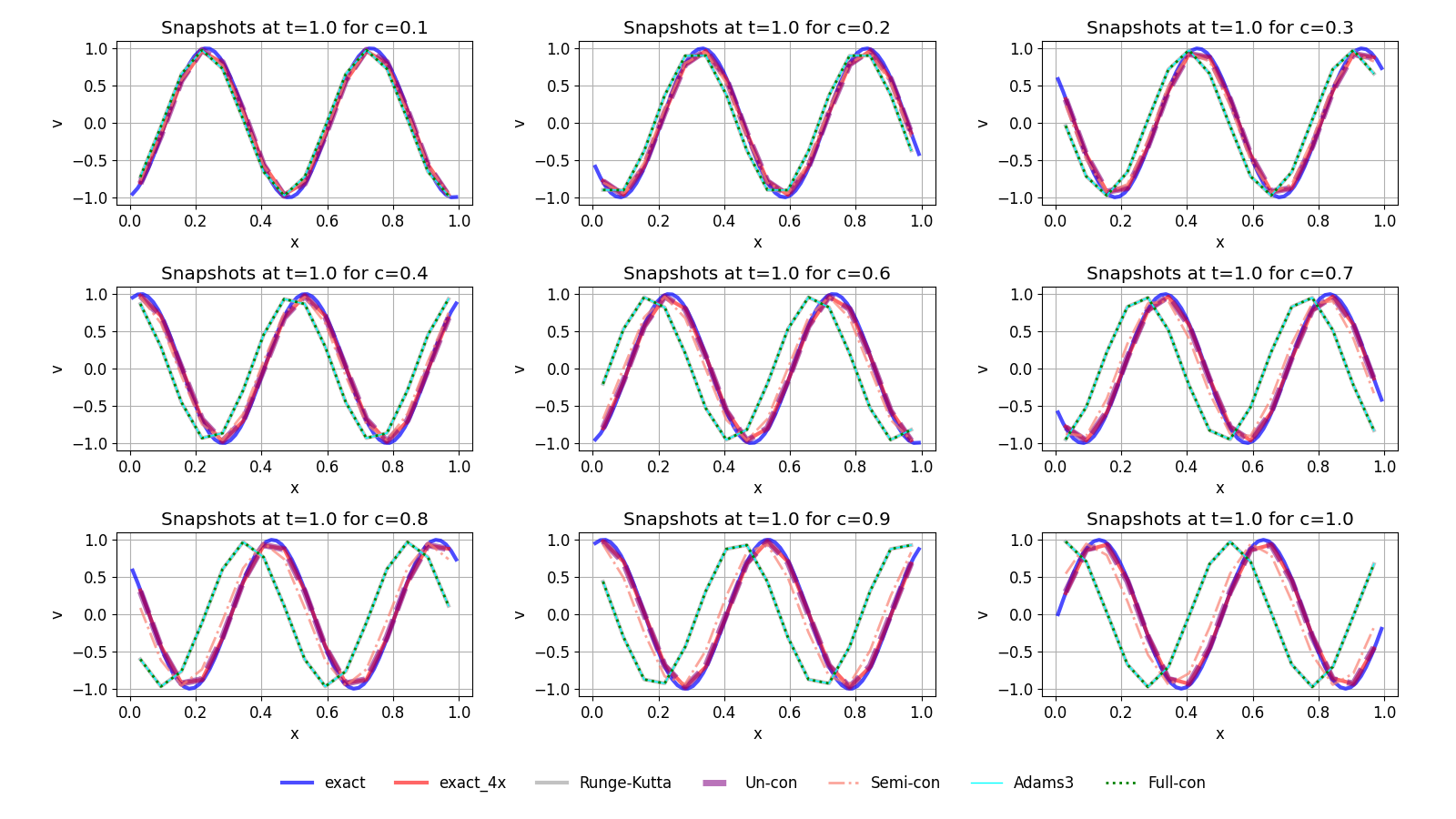}
\caption{The snapshots of wave equation for $c \in \left \{ 0.1, 0.2, 0.3, 0.4, 0.6, 0.7, 0.8, 0.9, 1.0 \right \} $ at $t=1$. The horizontal axis represents the spatial position, and the vertical axis represents the coarse-grained waves predicted by different methods at $t=1$.
}
\label{fig:wave_at1}
\end{figure}

In the following, we further analyze the error of the 1-D wave equation to provide some understanding of the improvement obtained from the learned time integration scheme. 
As shown from the Fourier analysis, the spatial discretization scheme introduces significant dispersion errors on coarse grids for the wave equation, here manifesting as a gradual lagging behind the phase during propagation. 
Figure \ref{fig:wave_at1} illustrates the solution of different models corresponding to different $c$ at $t=1$. It is evident that the phase error increases as $c$ increases. Among them, the Runge-Kutta method, the third-order Adams method, and the fully-constrained model lag in phase the most. 

Specifically, we aim to investigate how the coefficients of the unconstrained and semi-constrained models lead to error reductions. The third-order explicit linear multistep method combined with a second-order finite-volume method corresponds to a four-level explicit scheme. Assuming that the solution $v^{n-2}, v^{n-1}$ and $v^{n}$ from the three previous time steps are given by the exact solution Eq.~(\ref{eq:wave_exact}), we intend to examine the phase error brought about by the discretization scheme after advancing for one time step to $n+1$. 

\begin{figure}[htb!]
\centering
\includegraphics[width=14cm]{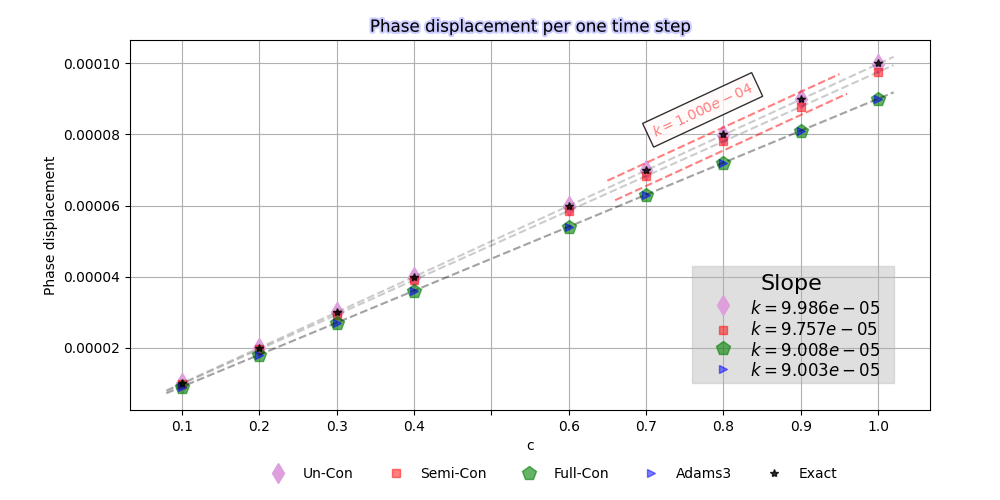}
\caption{Phase displacement per one time step for different time integration schemes for 1-D wave equation with different values of $c$. The grey dashed lines from top to bottom are the lines fitted by the unconstrained, semi-constrained, fully-constrained, and third-order Adams methods, with the latter two almost overlapping. Two pink dashed lines with $k = 1.000e-4$ is the slope for the exact solution, which equals to the size of the time step. The values in the grey box correspond to the slopes of various lines in the figure. Since $c=0.5$ is used to generate the training data, the corresponding results are not included.
}
\label{fig:mean_phase}
\end{figure}

It can be proven that the phase displacement per one time step is a constant for the coefficients given by the linear multistep method {\color{black}for initial conditions in Eq~(\ref{eq:wave_eq})}. Figure \ref{fig:mean_phase} shows the phase displacement per one time step for cases with different $c$. Since the coefficients from the data-driven models vary with the inputs, we have plotted the time-averaged phase displacements for these models. Ideally, the phase would move $c\Delta t$ in one time step, with the exact slope being of $\Delta t$ as shown in Figure \ref{fig:mean_phase}. 
It is seen that phase displacement from the unconstrained model is the closest to the exact one for different values of $c$. 
Improvements are also observed for the semi-constrained model. As for the fully-constrained model and third-order Adams method, however, phase displacements lower by approximately 10 percent are observed, being consistent with the observations in Figure \ref{fig:wave_at1}. Overall, Figure \ref{fig:mean_phase} demonstrates that the learned time integration scheme is capable of correcting the numerical phase displacement, which reduces the dispersion error caused by the spatial discretizations. More details on the mathematical derivations can be found in Appendix \ref{sec:Phase}.

\subsection{1-D Burgers' Equation} \label{sec:burgers}
In this section, we present the results of Burgers' equation from the data-driven time integration schemes and compare them with those from the Runge-Kutta method of order 3(2). Three different test cases, including different forcing terms, long integration time and large computational domain are considered. The performance of the data-driven time integration schemes is evaluated by applying them to simulations on a coarse grid, which is 32$\times$ coarser than the fine grid simulations (512 grids points) employed for generating the training data. In the 32$\times$ coarse grid simulations, the data-driven spatial disretization schemes trained (the number of steps and the learning rate for training are slightly different from the reference) using the method in the reference~\cite{2018Learning} are employed. 

The considered one-dimensional Burgers' equation is in the following form,
\begin{linenomath*}
\begin{align}
    \frac{\partial v}{\partial t}+ \frac{\partial}{\partial x} \left(\frac{v^{2} }{2}-\eta\frac{\partial v}{\partial x} \right)=f(x,t),
    \label{eq:burgers_eq}
\end{align}
\end{linenomath*}
 where $\eta=0.01$ is the viscosity and $f$ is the forcing term. The initial condition is $v(x,t=0)=0$. The periodic boundary condition is applied in $x$-direction. The forcing term $f$ is given as follows~\cite{2018Learning},
\begin{linenomath*}
\begin{equation}
    f(x,t) = \sum_{i=1}^{20}A_i \sin (\omega_i t + 2\pi l_i x/L + \phi_i),
\end{equation}
\end{linenomath*}
where $A_{i}\in [-0.5, 0.5]$, $\omega_{i} \in [-0.4, 0.4]$, $\phi_{i} \in [0,2\pi]$ and $l_{i}\in \{3,4,5,6\}$. The initial weights are uniformly distributed in the range of $-0.0001$ to $0.0001$ and $-0.001$ to $0.001$ for the unconstrained model and the semi-constrained mode respectively, with the same initial biases $[0, 0, -1, 1/2, -4/3, 23/12]$. The initial biases of the fully-constrained model are set as $[0, 0, 1/3, -4/3]$, and the initial weights are uniformly distributed in the range of $-0.01$ to $0.01$.  The $\gamma$ in loss function Eq.~(\ref{eq:loss}) is set to $10^{-15}$ for the semi-constrained model and $10^{-12}$ for the fully-constrained model. The three coarse-grained solutions of the Burgers' equation for $t\in \left [ 0, 20 \right ]$ with different forcing terms in the training set do not intersect with the test cases.

\begin{figure}[htb!]
\centering
\includegraphics[width=14cm]{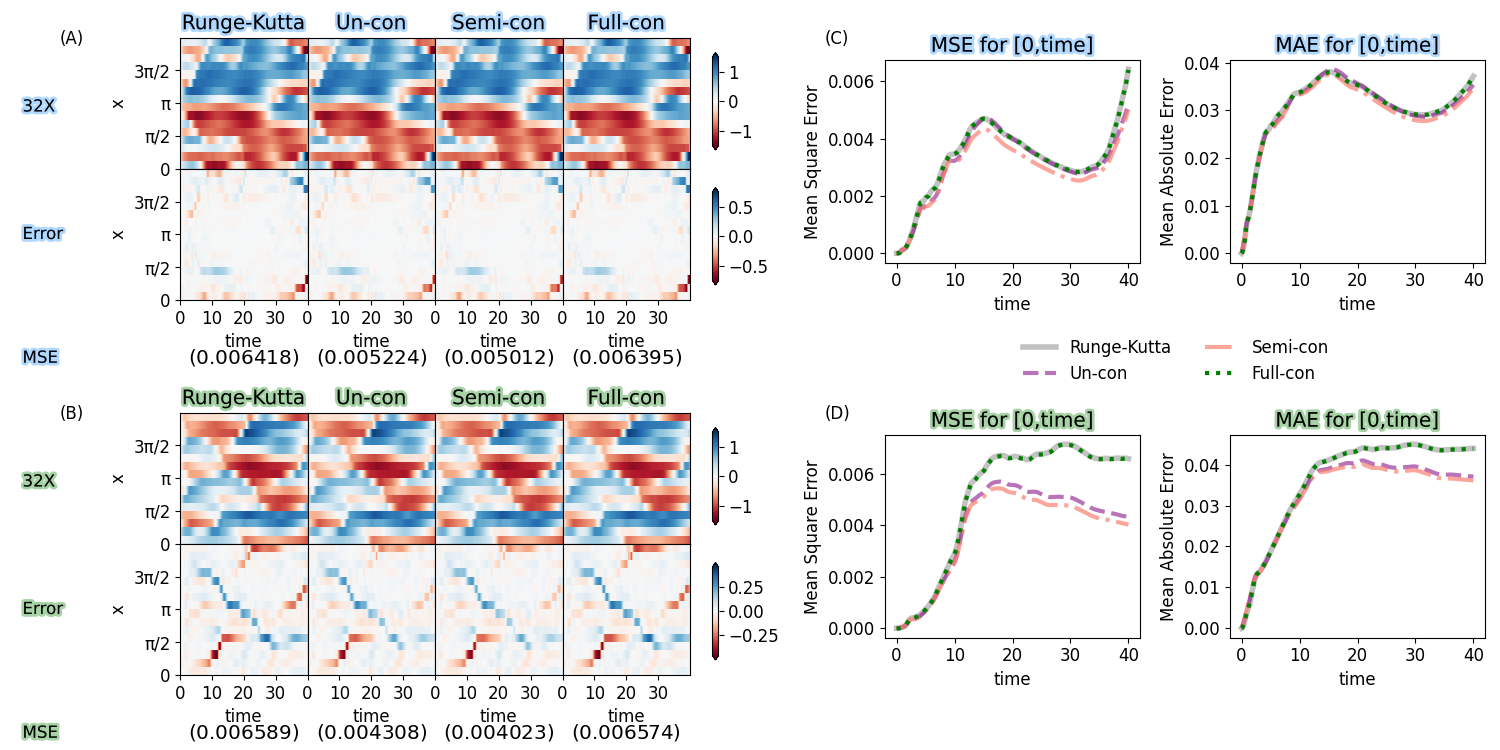}
	  \caption{Test results for two distinct forcing terms for (A, B) the realizations of solutions and error distribution and (C, D) the corresponding mean square error curves and mean absolute error curves. The employed grid is 32 times coarser than the reference fine grid. 
   The numbers in brackets below each subgraph in (A, B) are the mean square error averaged over the whole domain. The error shown in (C, D) is obtained by averaging the error in space ($[0, 2\pi]$) and time ($[0,t]$)
      }\label{fig:40_2pi}
\end{figure}

First, we examine the performance of the three data-driven time integration schemes using the test cases with the forcing term different from the training dataset. Figure \ref{fig:40_2pi} shows the results from two typical cases with (B, D with green-background title) and without (A, C with blue-background title) significant improvements, respectively. It is seen that the predictions from different time integration schemes are basically the same in the beginning of time. As further advancing in time, the unconstrained model and the semi-constrained model perform better than the other two models, which is measured using MSE and MAE. Little difference is observed between the predictions from the fully-constrained model and Runge-Kutta method, as no much space is left for further improvement when the constraints for deriving a third-order method are fully enforced when training the data-driven schemes. 

\begin{figure}[htb!]
\centering
\includegraphics[width=14cm]{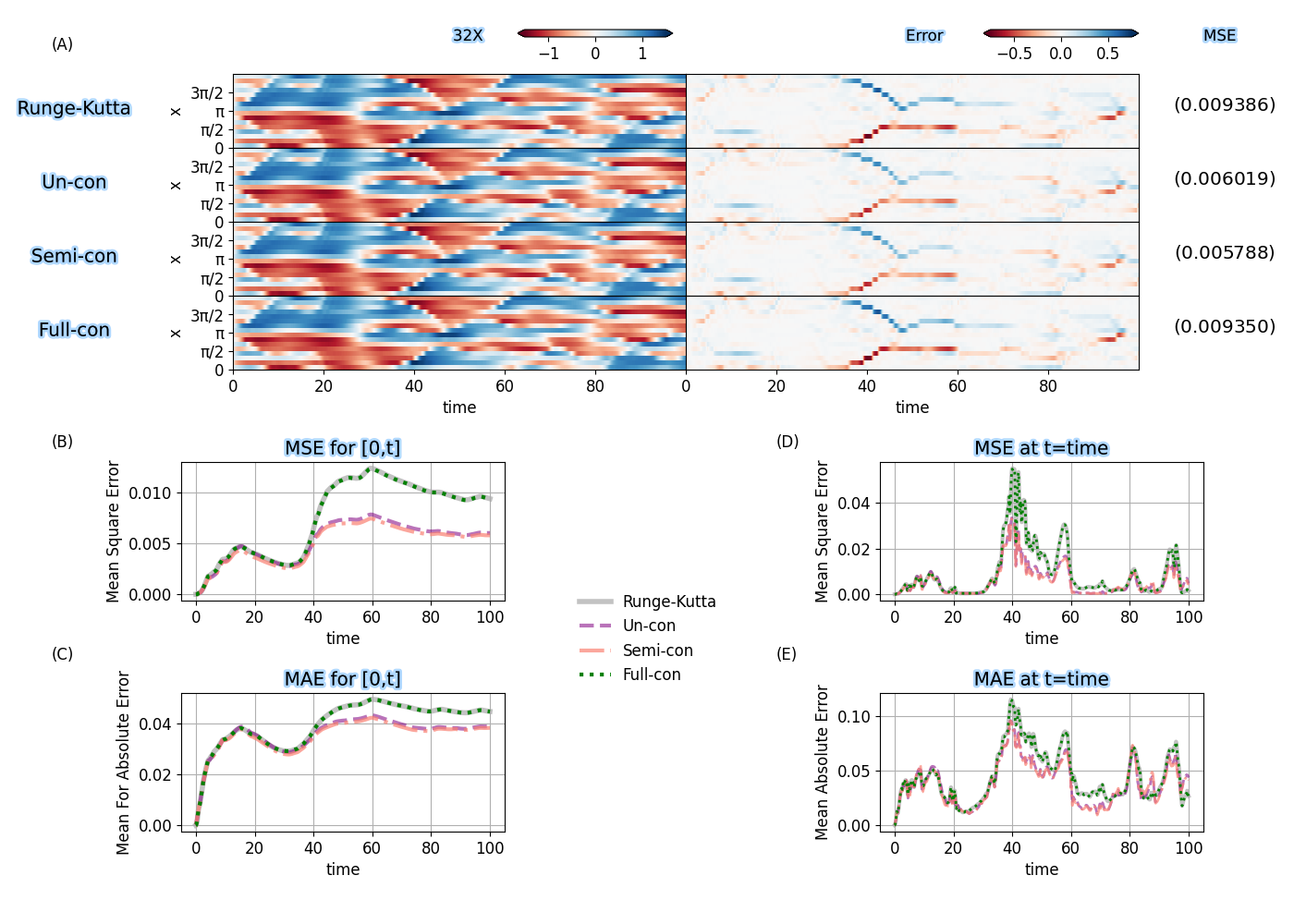}
	  \caption{Test results for a long integration time with $t\in \left [ 0,100 \right ] $ for (A) a realization of the solution and the corresponding error distribution in space and time, (B,D) the mean square error curves, and (C,E) the mean absolute error curves. The employed grid is 32 times coarser than {\color{black}the reference fine grid.} 
   The error shown in (B, C) is obtained by averaging the error in space ($[0, 2\pi]$) and time ($[0,t]$). The error shown in (D,E) is the error averaged over space ($[0, 2\pi]$) at $t$ instant.
      }\label{fig:100_2pi}
\end{figure}

\begin{figure}[htb!]
\centering
\includegraphics[width=14cm]{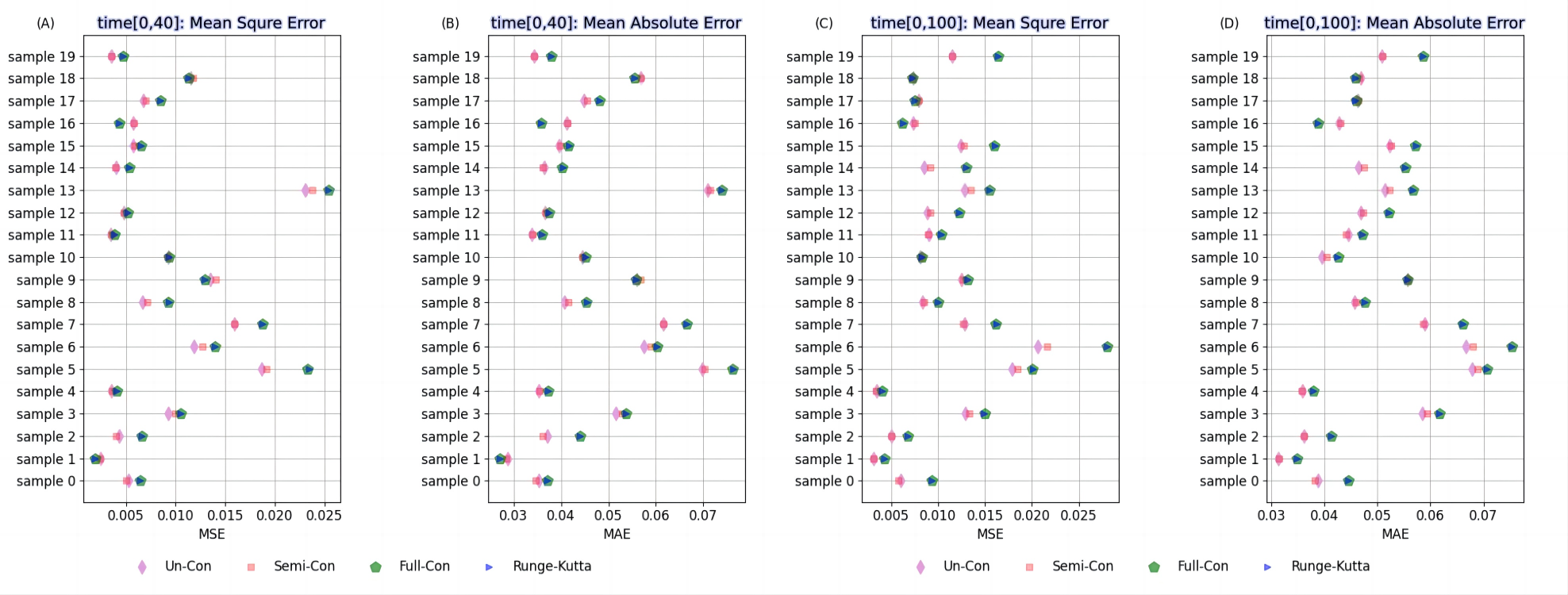}
\caption{Mean square error and mean absolute error of 20 samples for comparison with the time period of 40 and 100. (A) is the mean square error of these samples ($ 0\le t\le 40$) solved by four different time integration schemes. (B) is the mean absolute error of samples ($ 0\le t\le 40$) solved by four different time integration schemes. (C) and (D) are the same as (A) and (B), but for $ 0\le t\le 100$.}
\label{fig:20eq}
\end{figure}

A greater expectation of the proposed data-driven time integration scheme is their superior performance over an integration time longer than that of the training dataset. Figure \ref{fig:100_2pi} compares the predictions from the three learned time integration schemes trained on a temporal domain with $t\in \left [ 0,20 \right ] $ with those from the explicit Runge-Kutta method on a temporal domain with $t\in \left [ 0,100 \right ] $. It is clear that the unconstrained model and semi-constrained model outperform the fully-constrained model and the Runge-Kutta method as shown by the MSE and MAE for a long integration time.
In order to systematically evaluate the performance of the data-driven time integration schemes, 20 cases with vastly different random forcing terms are carried out. The exact value of error and the percentage of the error reduction are shown in the Appendix\ref{sec:Burger' Details }. As shown in Figure \ref{fig:20eq}, the error of learned unconstrained (purple diamond) and semi-constrained (orange square) discretization schemes, which are close to each other, are less than that of the Runge-Kutta method (blue triangle) in terms of MSE and MAE for most cases. As for the fully-constrained discretization scheme, the values of the MSE and MAE (green pentagon) are approximately the same with those from the Runge-Kutta method. 

Furthermore, the differences between the error from the three data-driven time integration schemes and the errors from the Runge-Kutta method are analysed by the paired t-tests. The P-values are shown in table \ref{tab:p-value} for the three schemes. At a significance level of $0.01$, it is seen that the error from the learned unconstrained and semi-constrained time integration schemes is significantly different from (lower than) the error from the Runge-Kutta method.

\begin{table}[htb!]
\centering
\caption{{\color{black}The average of mean square errors and the P-values of paired t-test of the mean square error and the mean absolute error of three data-driven time integration schemes} (The mean MSE is $0.011989284$ for the Runge-Kutta method).}

\label{tab:p-value}
\begin{tabular}{cccc} 
\toprule
            & Un-Con      & Semi-Con    & Full-Con     \\ 
\midrule
Mean MSEs   & 0.009832659 & 0.010041047 & 0.011981383  \\
MSE p-value & 0.000147995 & 0.000201141 & 0.087204851  \\
MAE p-value & 0.000111214 & 0.000180005 & 0.008416531  \\
\bottomrule
\end{tabular}
\end{table}

Overall, the test results from these cases with different forcing terms and a long integration time ($t\in [0,100]$ compared with $t\in [0,20]$ for the training dataset) demonstrate that the data-driven unconstrained and semi-constrained time integration schemes can reduce the error caused by the spatial coarse-graining. However, the fully-constrained model can hardly reduce the error, for which we deduce that fully enforcing the constraint conditions, which leaves very little room for improvement when training the model, is the key reason. 

So far, the size of the spatial domain ($[0,2\pi]$) employed in the test cases is the same as that of the training dataset. In the following, we test the data-driven time integration schemes using test cases with a larger spatial domain, in which it will use part of the spatial points as input to determine the coefficients of time integration schemes. Specifically, a 10$\times$ domain is employed in these test cases, that the periodic boundary condition is applied for $x \in \left [ 0,20\pi  \right ]$ and the forcing term $f$ should be  modified as reference~\cite{2018Learning}. This poses a challenge on selecting the grid points where the solutions are employed as the input of the data-driven model. In the tests of this work, the numerical solution on the first 16 grid points for $x \in \left [ 0,2\pi  \right ]$ are fed into the data-driven time integration schemes. The length of the integration time is 40. 

\begin{figure}[htb!]
\centering
\includegraphics[width=14cm]{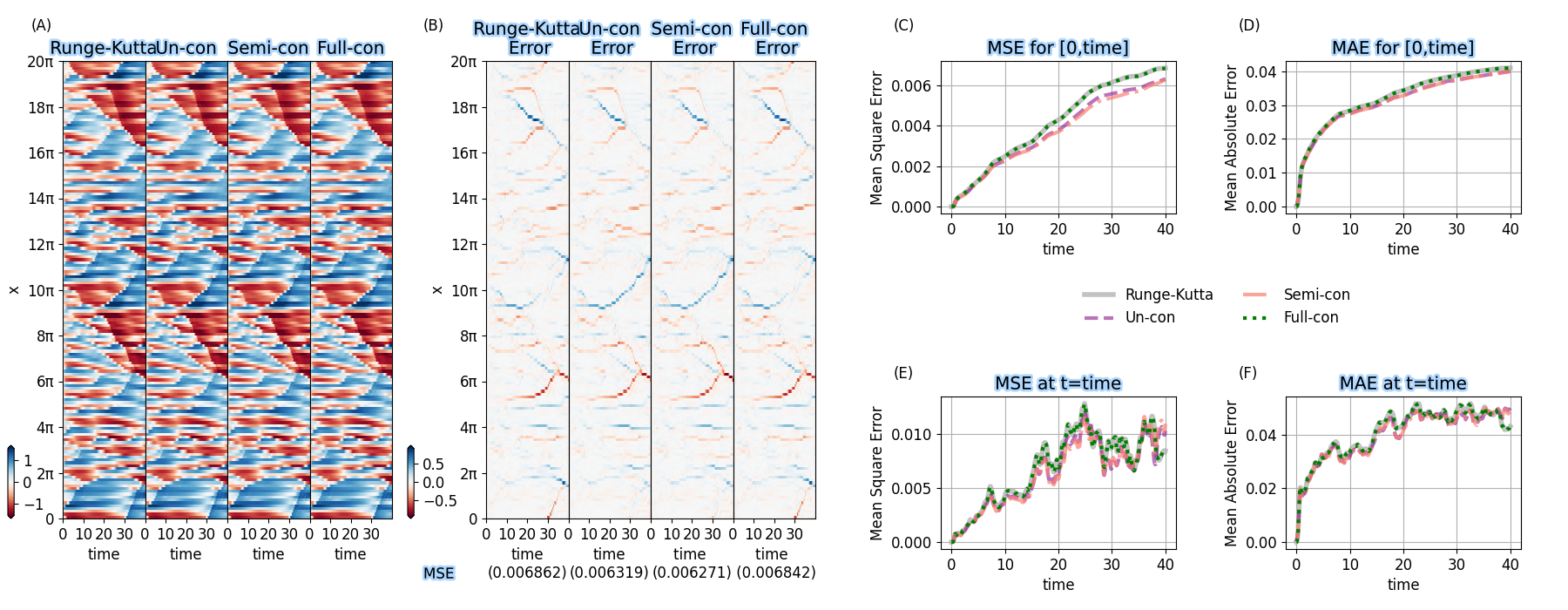}
\caption{Test results on a 10$\times$ larger spatial domain for (A) a realization of the solution, (B) the corresponding error distribution in space and time, (C,E) the mean square error curves, and (D,F) the mean absolute error curves. The spatial resolution is the same as the $32\times$ coarse test for spatial domain $[0, 20\pi]$. The error shown in (C, D) is obtained by averaging the error in space ($[0, 20\pi]$) and time ($[0,t]$).The error shown in (E,F) is the error averaged over space ($[0, 20\pi]$) at $t$ instant.
}
\label{fig:20pi}
\end{figure}

\begin{figure}[htb!]
\centering
\includegraphics[width=14cm]{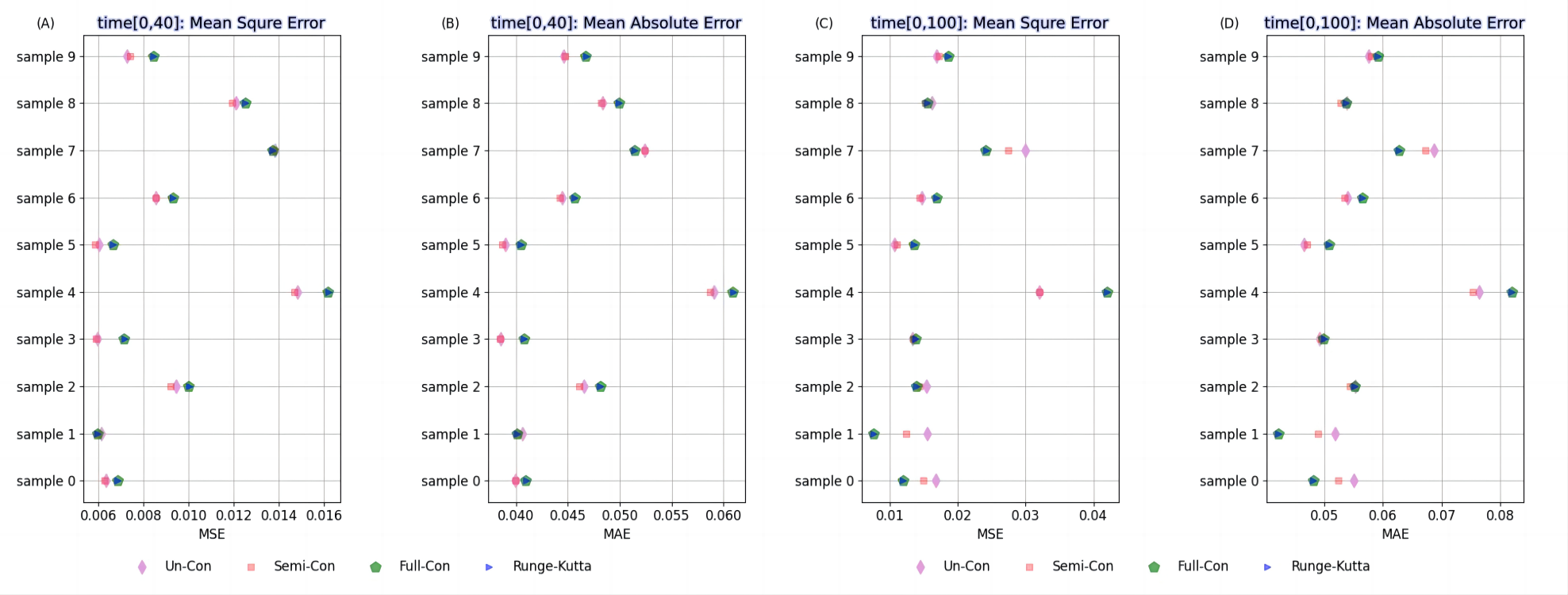}
\caption{Mean square error and mean absolute error of 10 samples for comparison on a 10$\times$ spatial domain with the time period for (A, B) $t\in[0,40]$ and (C, D) $t\in[0, 100]$ .}
\label{fig:20pi_10eqs}
\end{figure}

Figure \ref{fig:20pi} illustrates one experiment of the prediction on $x \in \left [ 0,20\pi  \right ]$. It is seen that the unconstrained and the semi-constrained models still outperform the Runge-Kutta method and the fully-constrained model, keeping lower values of MSE and MAE as advancing in time. Figure \ref{fig:20pi_10eqs} shows the error of 10 tests in the domain of $\left [ 0,20\pi  \right ]\times\left [ 0 ,40\right ]$ and $\left [ 0,20\pi  \right ]\times\left [ 0 ,100\right ]$, respectively. As seen, an overall improvement is obtained, although the {\color{black}performance deteriorates} for several cases as the time advances further to 100. 

Given the small initial weights of our neural network, the time-varying coefficients are approximately constant after retaining several decimal places. For data-driven unconstrained and semi-constrained time integration schemes, the constant optimized coefficients are used to calculate the above 20 samples on the domain $t\in \left [ 0,100 \right ]$, $x\in \left [ 0,2\pi \right ] $. Paired t-tests similar to the above are carried out. The exact value of error are shown in the Appendix\ref{sec:Burger' Details } and the p-values are shown in table \ref{tab:p-value(constant)}. {\color{black}It is seen that the constant optimized coefficients can produce results similar with the time-varying coefficients}. Although the p-values are slightly changed, it still can be considered acceptable at the significance level of 0.01. From the perspective of computing efficiency, the obtained constant optimized coefficients can be employed.

\begin{table}[htb!]
\centering
\caption{The average of mean square errors and the P-values of paired t-test of the mean square error and the mean absolute error of the constant optimized coefficients obtained from the unconstrained and semi-constrained models.}
\label{tab:p-value(constant)}
\begin{tabular}{ccc} 
\toprule
           & Const-coef (Un-Con) & Const-coef (Semi-Con)  \\ 
\midrule
Mean MSEs   & 0.009864970  & 0.009423226 \\
MSE p-value & 0.000160390  & 0.000341949 \\
MAE p-value & 0.000120366  & 0.000350297 \\
\bottomrule
\end{tabular}
\end{table}

Last but not least, we attempt to explore the reason for the improved performance of the data-driven time integration schemes. Can the improvements be obtained by increasing the scheme's order of accuracy? To probe into this, the Adams-Bashforth schemes of different orders of accuracy as shown in Table \ref{tab:Adams-Bashforth} are tested and compared with the constant optimized coefficients obtained from the unconstrained model. 

\begin{figure}[htb!]
\centering
\includegraphics[width=14cm]{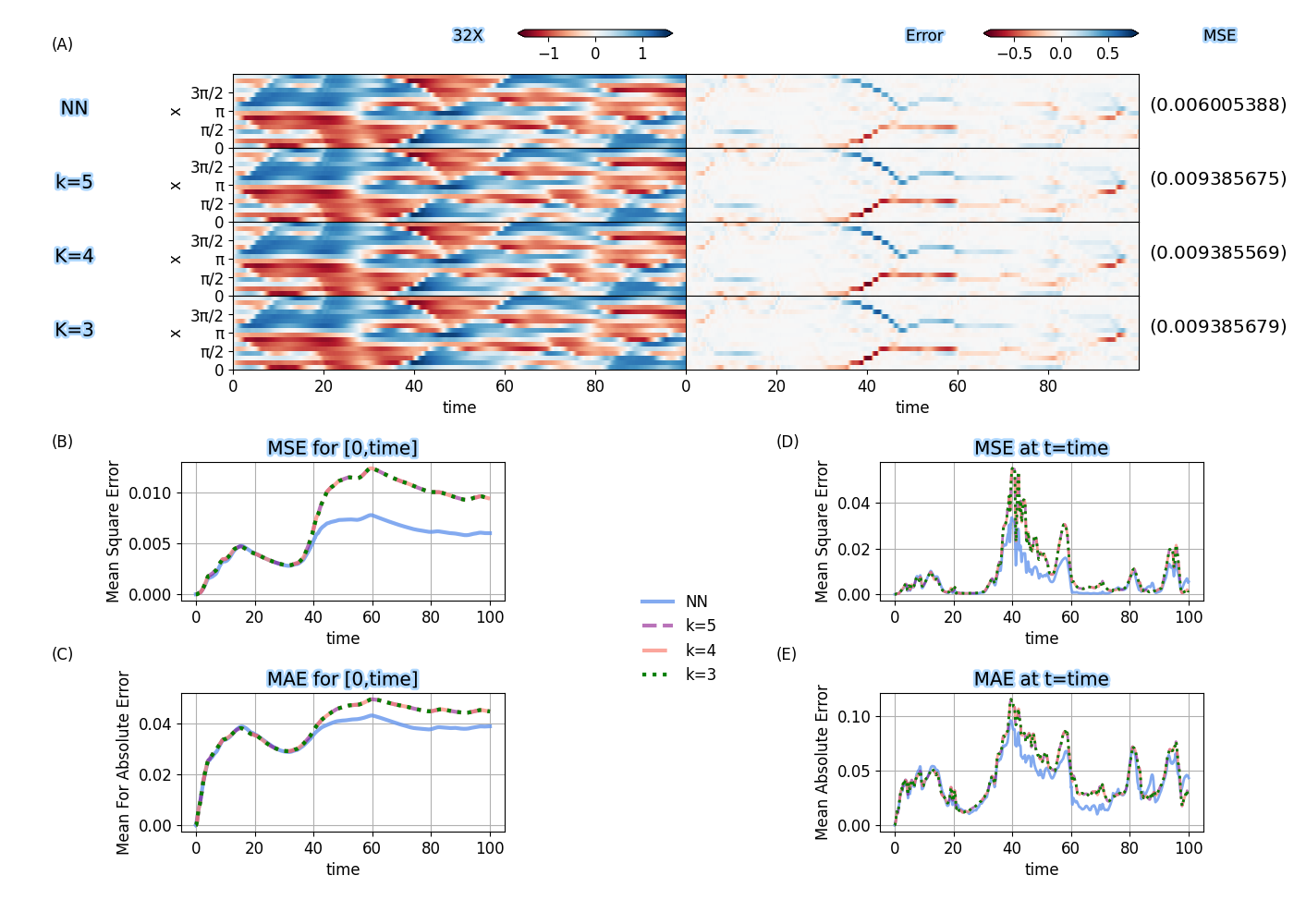}
	  \caption{Test results of Adams-Bashforth (k=3,4,5) and constant optimized coefficients  (Un-Con) for 
     (A) a realization of the solution and the corresponding error distribution in space and time, (B,D) the mean square error curves, and (C,E) the mean absolute error curves. The meaning of subgraphs are the same as Figure \ref{fig:100_2pi}.
     }\label{fig:Adams345}
\end{figure}

As shown in the figure \ref{fig:Adams345} for the results from a sample, improving the order of accuracy of the time integration schemes does not reduce the error during the advancement in time. The low spatial resolution employed in coarse grained simulations is the major source for error. The data-driven time integration scheme provides a new mechanism for canceling the error due to the coarse graining in space, at the price of destroying some properties or/and conditions employed for developing the conventional time integration schemes. Therefore, it is not surprising that imposing less or no constraints during model training, which leaves more space for optimizing the scheme coefficients, achieves an overall better performance. 
However, it is likely that the data-driven time integration schemes learned for certain types of cases (e.g., one specific coarsening or one specific PDE) may not improve the performance for other cases, as the learned schemes are uniquely tuned to reduce the error of a specific grid resolution for a specific PDE discretized in space using a specific scheme. 
The worst case is that the data-driven schemes have to be learned case by case, which is difficult for problems requiring a large amount of computational resources, e.g., high Reynolds number turbulent flows. This issue of generalization ability exists for almost all data-driven models. Further investigations need to be carried out in future work. 

\begin{table}[htb!]
\centering
\caption{Adams-Bashforth schemes for $k=3,4,5$ and the corresponding expression for truncation errors.}
\label{tab:Adams-Bashforth}
\begin{tabular}{ccc} 
\toprule
           k-step & Schemes      & Truncation error                                                                   \\ 
\midrule
$k=3$                                                & $v^{n+1}=v^n+\frac{\Delta t}{12}\left(23 F_{n}-16 F_{n-1}+5 F_{n-2}\right)$                               & $\frac{3}{8} (\Delta t)^4 \frac{\partial^{4} v}{\partial t^{4}} |_t=\zeta_n $      \\
$k=4$                                                & $v^{n+1}=v^n+\frac{\Delta t}{24}\left(55 F_{n}-59 F_{n-1}+37 F_{n-2}-9 F_{n-3}\right)$                      & $\frac{251}{720} (\Delta t)^5 \frac{\partial^{5} v}{\partial t^{5}} |_t=\zeta_n $  \\
$k=5$                                                & $v^{n+1}= v^n+\frac{\Delta t}{720}\left(1901 F_{n}-2774 F_{n-1}+2616 F_{n-2}-1274 F_{n-3}+251 F_{n-4}\right)$ & $\frac{95}{288} (\Delta t)^6 \frac{\partial^{6} v}{\partial t^{6}} |_t=\zeta_n $   \\
\bottomrule
\end{tabular}
\end{table}

\section{Conclusions}\label{sec:conclusions}
In this work, we proposed to learn time integration schemes using neural networks for solving partial differential equations on coarse grids, and tested the learned 3-step linear multistep method using the one-dimensional heat equation, the one-dimensional wave equation, and the one-dimensional Burgers' equation. During the training of the model, mathematical constraints, i.e., the consistency condition and the root condition (stability condition), are enforced. 
A backpropagation neural network with three layers and a low learning rate was employed with the initial values of the coefficients given as those of the conventional time integration schemes with perturbations. Three distinct time integration schemes were trained using the unconstrained model, the semi-constrained model with the root condition enforced, the fully-constrained model with both root and consistency conditions enforced.

The test results showed that the time integration schemes learned using the semi-constrained model and the unconstrained model are capable of reducing the mean square error (MSE) and the mean absolute error (MAE) for most cases, showing a reduction in error as high as an order of magnitude for the 1-D heat and the 1-D wave equation, and the most reduction in error in the range of 35\% to 40\% for the 1-D Burgers' equation. For the fully-constrained model, the prediction errors are close to those of the conventional time integration schemes. 

Analysis of the 1-D wave equation case revealed that the learned scheme effectively mitigates the dispersion error induced by the coarse grid. 
Further analysis of results from the Burgers' equation indicates that, instead of the order of accuracy, the data-driven model learned a mechanism to offset the error due to low spatial resolutions. 
Such mechanism is not readily analyzed using the existing methods such as Taylor analysis or Fourier analysis. 
To develop discretization schemes that exhibit strong generalization and interpretability, future research should focus on developing mathematical theories and tools for analyzing and evaluating these learned schemes.
 
\section*{Acknowledgment}\label{sec:acknowledgement}
This work was funded by the NSFC Basic Science Center Program for ``Multiscale Problems in Nonlinear Mechanics'' (NO. 11988102), National Natural Science Foundation of China (NO. 12172360), Institute of Mechanics CAS, and Chinese Academy of Sciences.

\linespread{0.5} \selectfont
\bibliographystyle{unsrt}
\bibliography{refs}

\begin{thebibliography}{10}

\bibitem{smagorinsky1963general}
Joseph Smagorinsky.
\newblock General circulation experiments with the primitive equations: I. the basic experiment.
\newblock {\em Monthly Weather Review}, 91(3):99--164, 1963.

\bibitem{germano1991dynamic}
Massimo Germano, Ugo Piomelli, Parviz Moin, and William~H Cabot.
\newblock A dynamic subgrid-scale eddy viscosity model.
\newblock {\em Physics of Fluids A: Fluid Dynamics}, 3(7):1760--1765, 1991.

\bibitem{grinstein2007implicit}
Fernando~F Grinstein, Len~G Margolin, and William~J Rider.
\newblock {\em Implicit large eddy simulation}, volume~10.
\newblock Cambridge University Press Cambridge, 2007.

\bibitem{Wyngaard2004}
JC~Wyngaard.
\newblock Toward numerical modeling in the ``terra incognita''.
\newblock {\em Journal of the Atmospheric Sciences}, 61(14):1816--1826, JUL 2004.

\bibitem{Honnert2020}
Rachel Honnert, Georgios~A. Efstathiou, Robert~J. Beare, Junshi Ito, Adrian Lock, Roel Neggers, Robert~S. Plant, Hyeyum~Hailey Shin, Lorenzo Tomassini, and Bowen Zhou.
\newblock The atmospheric boundary layer and the ``gray zone{''} of turbulence: A critical review.
\newblock {\em Journal of Geophysical Research-Atmospheres}, 125(13), JUL 16 2020.

\bibitem{Raissi2017a}
Maziar Raissi, Paris Perdikaris, and George~Em Karniadakis.
\newblock Physics informed deep learning (part i): Data-driven solutions of nonlinear partial differential equations.
\newblock {\em ar{X}iv:1711.10561}, 2017.

\bibitem{Sirignano2018}
Justin Sirignano and Konstantinos Spiliopoulos.
\newblock Dgm: A deep learning algorithm for solving partial differential equations.
\newblock {\em Journal of Computational Physics}, 375:1339--1364, DEC 15 2018.

\bibitem{Lu2021}
Lu~Lu, Xuhui Meng, Zhiping Mao, and George~Em Karniadakis.
\newblock Deepxde: A deep learning library for solving differential equations.
\newblock {\em SIAM Review}, 63(1):208--228, MAR 2021.

\bibitem{Wang2021}
Sifan Wang, Hanwen Wang, and Paris Perdikaris.
\newblock On the eigenvector bias of fourier feature networks: From regression to solving multi-scale pdes with physics-informed neural networks.
\newblock {\em Computer Methods in Applied Mechanics and Engineering}, 384, OCT 1 2021.

\bibitem{Zhang2019}
Dongkun Zhang, Ling Guo, and George~Em Karniadakis.
\newblock Learning in modal space: Solving time-dependent stochastic pdes using physics-informed neural networks.
\newblock {\em ar{X}iv:1905.01205v2}, 2019.

\bibitem{ren2022phycrnet}
Pu~Ren, Chengping Rao, Yang Liu, Jian-Xun Wang, and Hao Sun.
\newblock Phycrnet: Physics-informed convolutional-recurrent network for solving spatiotemporal pdes.
\newblock {\em Computer Methods in Applied Mechanics and Engineering}, 389:114399, 2022.

\bibitem{Karniadakis2021}
George~Em Karniadakis, Ioannis~G. Kevrekidis, Lu~Lu, Paris Perdikaris, Sifan Wang, and Liu Yang.
\newblock Physics-informed machine learning.
\newblock {\em Nature Reviews Physics}, 3(6):422--440, JUN 2021.

\bibitem{Cuomo2022}
Salvatore Cuomo, Vincenzo~Schiano Di~Cola, Fabio Giampaolo, Gianluigi Rozza, Maziar Raissi, and Francesco Piccialli.
\newblock Scientific machine learning through physics-informed neural networks: Where we are and what's next.
\newblock {\em Journal of Scientific Computing}, 92(3), SEP 2022.

\bibitem{Poggio2017}
Tomaso Poggio, Hrushikesh Mhaskar, Lorenzo Rosasco, Brando Miranda, and Qianli Liao.
\newblock Why and when can deep-but not shallow-networks avoid the curse of dimensionality: A review.
\newblock {\em International Journal of Automation and Computing}, 14(5, SI):503--519, OCT 2017.

\bibitem{Blechschmidt2021}
Jan Blechschmidt and Oliver~G. Ernst.
\newblock Three ways to solve partial differential equations with neural networks — a review.
\newblock {\em GAMM-Mitteilungen}, 44:e202100006, 2021.

\bibitem{Raissi2017b}
Maziar Raissi, Paris Perdikaris, and George~Em Karniadakis.
\newblock Physics informed deep learning (part ii): Data-driven discovery of nonlinear partial differential equations.
\newblock {\em ar{X}iv:1711.10566}, 2017.

\bibitem{Raissi2019}
Maziar Raissi, Paris Perdikaris, and George~Em Karniadakis.
\newblock Physics-informed neural networks: A deep learning framework for solving forward and inverse problems involving nonlinear partial differential equations.
\newblock {\em Journal of Computational Physics}, 378:686--707, 2019.

\bibitem{Beck2018}
Christian Beck, Sebastian Becker, Philipp Grohs, Nor Jaafari, and Arnulf Jentzen.
\newblock Solving the kolmogorov pde by means of deep learning.
\newblock {\em ar{X}iv:1806.00421}, 2018.

\bibitem{E2017}
Weinan E, Jiequn Han, and Arnulf Jentzen.
\newblock Deep learning-based numerical methods for high-dimensional parabolic partial differential equations and backward stochastic differential equations.
\newblock {\em Communications in Mathematics and Statistics}, 5(4):349--380, DEC 2017.

\bibitem{Han2018}
Jiequn Han, Arnulf Jentzen, and E.~Weinan.
\newblock Solving high-dimensional partial differential equations using deep learning.
\newblock {\em Proceedings of the National Academy of Sciences of the United States of America}, 115(34):8505--8510, AUG 21 2018.

\bibitem{duraisamy2019turbulence}
Karthik Duraisamy, Gianluca Iaccarino, and Heng Xiao.
\newblock Turbulence modeling in the age of data.
\newblock {\em Annual Review of Fluid Mechanics}, 51:357--377, 2019.

\bibitem{zhu2019machine}
Linyang Zhu, Weiwei Zhang, Jiaqing Kou, and Yilang Liu.
\newblock Machine learning methods for turbulence modeling in subsonic flows around airfoils.
\newblock {\em Physics of Fluids}, 31(1):015105, 2019.

\bibitem{zhou2019subgrid}
Zhideng Zhou, Guowei He, Shizhao Wang, and Guodong Jin.
\newblock Subgrid-scale model for large-eddy simulation of isotropic turbulent flows using an artificial neural network.
\newblock {\em Computers \& Fluids}, 195:104319, 2019.

\bibitem{meng2022artificial}
Qingjia Meng, Zhou Jiang, and Jianchun Wang.
\newblock Artificial neural network-based subgrid-scale models for les of compressible turbulent channel flow.
\newblock {\em Theoretical and Applied Mechanics Letters}, page 100399, 2022.

\bibitem{Yang_etal_PRF_2019}
X.~I.~A. Yang, S.~Zafar, J.-X. Wang, and H.~Xiao.
\newblock Predictive large-eddy-simulation wall modeling via physics-informed neural networks.
\newblock {\em Physical Review Fluids}, 4:034602, 2019.

\bibitem{zhou2021wall}
Zhideng Zhou, Guowei He, and Xiaolei Yang.
\newblock Wall model based on neural networks for les of turbulent flows over periodic hills.
\newblock {\em Physical Review Fluids}, 6(5):054610, 2021.

\bibitem{Zhang2022}
Xin-Lei Zhang, Heng Xiao, Xiaodong Luo, and Guowei He.
\newblock Ensemble kalman method for learning turbulence models from indirect observation data.
\newblock {\em Journal of Fluid Mechanics}, 949, SEP 29 2022.

\bibitem{Zhang2023}
Xin-Lei Zhang, Heng Xiao, Xiaodong Luo, and Guowei He.
\newblock Combining direct and indirect sparse data for learning generalizable turbulence models.
\newblock {\em Journal of Computational Physics}, 489, SEP 15 2023.

\bibitem{bezgin2021data}
Deniz~A Bezgin, Steffen~J Schmidt, and Nikolaus~A Adams.
\newblock A data-driven physics-informed finite-volume scheme for nonclassical undercompressive shocks.
\newblock {\em Journal of Computational Physics}, 437:110324, 2021.

\bibitem{kossaczka2022neural}
T~Kossaczka, M~Ehrhardt, and M~G{\"u}nther.
\newblock A neural network enhanced weighted essentially non-oscillatory method for nonlinear degenerate parabolic equations.
\newblock {\em Physics of Fluids}, 34(2):026604, 2022.

\bibitem{2018Learning}
Yohai Bar-Sinai, Stephan Hoyer, Jason Hickey, and Michael~P. Brenner.
\newblock Learning data-driven discretizations for partial differential equations.
\newblock {\em Proceedings of the National Academy of Sciences}, 116(31):15344--15349, 2019.

\bibitem{hairer1993solving}
Ernst Hairer, Syvert~P N{\o}rsett, and Gerhard Wanner.
\newblock {\em Solving ordinary differential equations. 1, Nonstiff problems}.
\newblock Springer-Verlag, 1993.

\bibitem{gantmacher2005applications}
Feliks~Rouminovich Gantmacher and Joel~Lee Brenner.
\newblock {\em Applications of the theory of matrices}.
\newblock Courier Corporation, 2005.

\bibitem{luenberger1984linear}
David~G Luenberger, Yinyu Ye, et~al.
\newblock {\em Linear and nonlinear programming}, volume~2.
\newblock Springer, 1984.

\bibitem{anderson2016computational}
Dale Anderson, John~C Tannehill, and Richard~H Pletcher.
\newblock {\em Computational fluid mechanics and heat transfer}.
\newblock Taylor \& Francis, 2016.

\bibitem{lomax2001fundamentals}
Harvard Lomax, Thomas~H Pulliam, David~W Zingg, Thomas~H Pulliam, and David~W Zingg.
\newblock {\em Fundamentals of computational fluid dynamics}, volume 246.
\newblock Springer, 2001.

\end{thebibliography}

\begin{appendices}
\counterwithin{figure}{section}
\section{Detail results for heat equation}\label{app:heat}
In this appendix, the mean square error (MSE) and the mean absolute error (MAE) for the 1-D heat equation cases with different $\lambda$ are given in table~\ref{tab:heat_error}.

\begin{longtblr}[
  caption = {Results of 1-D heat equations for different $\lambda$. The MSE and MAE in the first four columns are the error averaged over space $([0, 1])$ and time $([0, 0.5])$ or $([0, 1])$. The max MSE and the max MAE in the last two columns are the max of the averaged error curve in space $([0, 1])$ and time $([0, t])$.},
  label={tab:heat_error},
]{
  cells = {c},
  cell{1}{1} = {c=7}{},
  cell{7}{1} = {c=7}{},
  cell{13}{1} = {c=7}{},
  cell{19}{1} = {c=7}{},
  cell{25}{1} = {c=7}{},
  cell{31}{1} = {c=7}{},
  cell{37}{1} = {c=7}{},
  cell{43}{1} = {c=7}{},
  cell{49}{1} = {c=7}{},
  vline{2} = {2-6,8-12,14-18,20-24,26-30,32-36,38-42,44-48,50-54}{},
  hline{1-3,7-9,13-15,19-21,25-27,31-33,37-39,43-45,49-51,55} = {-}{},
}

$\lambda$ = 0.1 &                  &                &                  &                &          &          \\
               & MSE for [0, 0.5] & MSE for [0, 1] & MAE for [0, 0.5] & MAE for [0, 1] & Max MSE  & Max MAE  \\
RK             & 7.82e-06         & 5.13e-06       & 2.43e-03         & 1.88e-03       & 7.94e-06 & 2.44e-03 \\
Un-con         & 2.12e-04         & 1.35e-04       & 1.27e-02         & 9.58e-03       & 2.17e-04 & 1.28e-02 \\
Semi-con       & 1.71e-05         & 1.11e-05       & 3.60e-03         & 2.76e-03       & 1.74e-05 & 3.62e-03 \\
Full-con       & 7.82e-06         & 5.13e-06       & 2.43e-03         & 1.88e-03       & 7.94e-06 & 2.44e-03 \\
$\lambda$ = 0.2 &                  &                &                  &                &          &          \\
               & MSE for [0, 0.5] & MSE for [0, 1] & MAE for [0, 0.5] & MAE for [0, 1] & Max MSE  & Max MAE  \\
RK             & 5.12e-06         & 2.60e-06       & 1.88e-03         & 1.04e-03       & 7.92e-06 & 2.44e-03 \\
Un-con         & 2.02e-05         & 1.02e-05       & 3.72e-03         & 2.04e-03       & 3.17e-05 & 4.88e-03 \\
Semi-con       & 5.48e-07         & 2.78e-07       & 6.14e-04         & 3.38e-04       & 8.52e-07 & 7.99e-04 \\
Full-con       & 5.05e-06         & 2.57e-06       & 1.87e-03         & 1.03e-03       & 7.82e-06 & 2.42e-03 \\
$\lambda$ = 0.3 &                  &                &                  &                &          &          \\
               & MSE for [0, 0.5] & MSE for [0, 1] & MAE for [0, 0.5] & MAE for [0, 1] & Max MSE  & Max MAE  \\
RK             & 3.44e-06         & 1.72e-06       & 1.36e-03         & 6.92e-04       & 7.85e-06 & 2.43e-03 \\
Un-con         & 2.75e-06         & 1.37e-06       & 1.21e-03         & 6.15e-04       & 6.33e-06 & 2.18e-03 \\
Semi-con       & 1.25e-08         & 6.26e-09       & 8.17e-05         & 4.16e-05       & 2.94e-08 & 1.50e-04 \\
Full-con       & 3.38e-06         & 1.69e-06       & 1.35e-03         & 6.86e-04       & 7.72e-06 & 2.41e-03 \\
$\lambda$ = 0.4 &                  &                &                  &                &          &          \\
               & MSE for [0, 0.5] & MSE for [0, 1] & MAE for [0, 0.5] & MAE for [0, 1] & Max MSE  & Max MAE  \\
RK             & 2.56e-06         & 1.28e-06       & 1.03e-03         & 5.17e-04       & 7.76e-06 & 2.41e-03 \\
Un-con         & 2.91e-07         & 1.46e-07       & 3.46e-04         & 1.74e-04       & 8.88e-07 & 8.14e-04 \\
Semi-con       & 1.69e-07         & 8.44e-08       & 2.64e-04         & 1.32e-04       & 5.20e-07 & 6.26e-04 \\
Full-con       & 2.51e-06         & 1.25e-06       & 1.02e-03         & 5.12e-04       & 7.64e-06 & 2.39e-03 \\
$\lambda$ = 0.6 &                  &                &                  &                &          &          \\
               & MSE for [0, 0.5] & MSE for [0, 1] & MAE for [0, 0.5] & MAE for [0, 1] & Max MSE  & Max MAE  \\
RK             & 1.64e-06         & 8.20e-07       & 6.78e-04         & 3.39e-04       & 7.40e-06 & 2.36e-03 \\
Un-con         & 8.92e-08         & 4.46e-08       & 1.57e-04         & 7.86e-05       & 4.14e-07 & 5.59e-04 \\
Semi-con       & 3.52e-07         & 1.76e-07       & 3.12e-04         & 1.56e-04       & 1.62e-06 & 1.10e-03 \\
Full-con       & 1.64e-06         & 8.18e-07       & 6.75e-04         & 3.37e-04       & 7.47e-06 & 2.37e-03 \\
$\lambda$ = 0.7 &                  &                &                  &                &          &          \\
               & MSE for [0, 0.5] & MSE for [0, 1] & MAE for [0, 0.5] & MAE for [0, 1] & Max MSE  & Max MAE  \\
RK             & 1.33e-06         & 6.67e-07       & 5.68e-04         & 2.84e-04       & 6.98e-06 & 2.29e-03 \\
Un-con         & 2.24e-07         & 1.12e-07       & 2.31e-04         & 1.15e-04       & 1.20e-06 & 9.53e-04 \\
Semi-con       & 3.81e-07         & 1.91e-07       & 3.01e-04         & 1.51e-04       & 2.04e-06 & 1.24e-03 \\
Full-con       & 1.39e-06         & 6.94e-07       & 5.75e-04         & 2.88e-04       & 7.39e-06 & 2.36e-03 \\
$\lambda$ = 0.8 &                  &                &                  &                &          &          \\
               & MSE for [0, 0.5] & MSE for [0, 1] & MAE for [0, 0.5] & MAE for [0, 1] & Max MSE  & Max MAE  \\
RK             & 1.14e-06         & 5.69e-07       & 4.92e-04         & 2.46e-04       & 6.72e-06 & 2.25e-03 \\
Un-con         & 3.38e-07         & 1.69e-07       & 2.65e-04         & 1.33e-04       & 2.07e-06 & 1.25e-03 \\
Semi-con       & 3.91e-07         & 1.96e-07       & 2.85e-04         & 1.43e-04       & 2.40e-06 & 1.34e-03 \\
Full-con       & 1.20e-06         & 6.00e-07       & 5.01e-04         & 2.50e-04       & 7.31e-06 & 2.34e-03 \\
$\lambda$ = 0.9 &                  &                &                  &                &          &          \\
               & MSE for [0, 0.5] & MSE for [0, 1] & MAE for [0, 0.5] & MAE for [0, 1] & Max MSE  & Max MAE  \\
RK             & 9.50e-07         & 4.75e-07       & 4.26e-04         & 2.14e-04       & 6.22e-06 & 2.17e-03 \\
Un-con         & 4.22e-07         & 2.11e-07       & 2.80e-04         & 1.40e-04       & 2.91e-06 & 1.48e-03 \\
Semi-con       & 3.91e-07         & 1.95e-07       & 2.69e-04         & 1.34e-04       & 2.69e-06 & 1.42e-03 \\
Full-con       & 1.06e-06         & 5.28e-07       & 4.43e-04         & 2.21e-04       & 7.24e-06 & 2.33e-03 \\
$\lambda $= 1.0 &                  &                &                  &                &          &          \\
               & MSE for [0, 0.5] & MSE for [0, 1] & MAE for [0, 0.5] & MAE for [0, 1] & Max MSE  & Max MAE  \\
RK             & 8.30e-07         & 4.15e-07       & 3.79e-04         & 1.90e-04       & 6.01e-06 & 2.13e-03 \\
Un-con         & 4.81e-07         & 2.41e-07       & 2.83e-04         & 1.42e-04       & 3.68e-06 & 1.66e-03 \\
Semi-con       & 3.84e-07         & 1.92e-07       & 2.53e-04         & 1.26e-04       & 2.94e-06 & 1.49e-03 \\
Full-con       & 9.40e-07         & 4.70e-07       & 3.96e-04         & 1.98e-04       & 7.16e-06 & 2.32e-03 
\end{longtblr}

\section{Detail results for wave equation}\label{app:wave}
The mean square error (MSE) and the mean absolute error (MAE) for the 1-D wave equation cases with different $c$ are given in table~\ref{tab:wave_error}.

\begin{longtblr}[
  caption = {Results of 1-D wave equations for different $c$. The explanations of the headers are the same as table~\ref{tab:heat_error}.},
  label={tab:wave_error},
]{
  cells = {c},
  cell{1}{1} = {c=7}{},
  cell{8}{1} = {c=7}{},
  cell{15}{1} = {c=7}{},
  cell{22}{1} = {c=7}{},
  cell{29}{1} = {c=7}{},
  cell{36}{1} = {c=7}{},
  cell{43}{1} = {c=7}{},
  cell{50}{1} = {c=7}{},
  cell{57}{1} = {c=7}{},
  vline{2} = {2-7,9-14,16-21,23-28,30-35,37-42,44-49,51-56,58-63}{},
  hline{1-3,8-10,15-17,22-24,29-31,36-38,43-45,50-52,57-59,64} = {-}{},
}
$c$ = 0.1  &                  &                &                  &                &          &          \\
         & MSE for [0, 0.5] & MSE for [0, 1] & MAE for [0, 0.5] & MAE for [0, 1] & Max MSE  & Max MAE  \\
RK       & 6.23e-04         & 2.49e-03       & 1.93e-02         & 3.86e-02       & 2.49e-03 & 3.86e-02 \\
Adams3   & 6.23e-04         & 2.49e-03       & 1.93e-02         & 3.86e-02       & 2.49e-03 & 3.86e-02 \\
Uncon    & 1.30e-05         & 4.53e-05       & 2.79e-03         & 5.31e-03       & 4.53e-05 & 5.31e-03 \\
Semi-con & 4.51e-05         & 1.76e-04       & 5.31e-03         & 1.04e-02       & 1.76e-04 & 1.04e-02 \\
Full-con & 6.23e-04         & 2.49e-03       & 1.93e-02         & 3.86e-02       & 2.49e-03 & 3.86e-02 \\
$c$ = 0.2  &                  &                &                  &                &          &          \\
         & MSE for [0, 0.5] & MSE for [0, 1] & MAE for [0, 0.5] & MAE for [0, 1] & Max MSE  & Max MAE  \\
RK       & 2.49e-03         & 9.93e-03       & 3.86e-02         & 7.79e-02       & 9.93e-03 & 7.79e-02 \\
Adams3   & 2.49e-03         & 9.93e-03       & 3.86e-02         & 7.79e-02       & 9.93e-03 & 7.79e-02 \\
Uncon    & 1.18e-05         & 2.34e-05       & 2.73e-03         & 4.01e-03       & 2.34e-05 & 4.01e-03 \\
Semi-con & 1.56e-04         & 6.15e-04       & 9.76e-03         & 1.94e-02       & 6.15e-04 & 1.94e-02 \\
Full-con & 2.49e-03         & 9.93e-03       & 3.86e-02         & 7.79e-02       & 9.93e-03 & 7.79e-02 \\
$c$ = 0.3  &                  &                &                  &                &          &          \\
         & MSE for [0, 0.5] & MSE for [0, 1] & MAE for [0, 0.5] & MAE for [0, 1] & Max MSE  & Max MAE  \\
RK       & 5.59e-03         & 2.23e-02       & 5.85e-02         & 1.16e-01       & 2.23e-02 & 1.16e-01 \\
Adams3   & 5.59e-03         & 2.23e-02       & 5.85e-02         & 1.16e-01       & 2.23e-02 & 1.16e-01 \\
Uncon    & 9.73e-06         & 1.23e-05       & 2.52e-03         & 2.98e-03       & 1.23e-05 & 2.98e-03 \\
Semi-con & 3.40e-04         & 1.35e-03       & 1.45e-02         & 2.86e-02       & 1.35e-03 & 2.86e-02 \\
Full-con & 5.59e-03         & 2.23e-02       & 5.85e-02         & 1.16e-01       & 2.23e-02 & 1.16e-01 \\
$c$ = 0.4  &                  &                &                  &                &          &          \\
         & MSE for [0, 0.5] & MSE for [0, 1] & MAE for [0, 0.5] & MAE for [0, 1] & Max MSE  & Max MAE  \\
RK       & 9.93e-03         & 3.94e-02       & 7.79e-02         & 1.55e-01       & 3.94e-02 & 1.55e-01 \\
Adams3   & 9.93e-03         & 3.94e-02       & 7.79e-02         & 1.55e-01       & 3.94e-02 & 1.55e-01 \\
Uncon    & 7.93e-06         & 1.23e-05       & 2.33e-03         & 2.98e-03       & 1.23e-05 & 2.98e-03 \\
Semi-con & 6.00e-04         & 2.39e-03       & 1.92e-02         & 3.82e-02       & 2.39e-03 & 3.82e-02 \\
Full-con & 9.93e-03         & 3.94e-02       & 7.79e-02         & 1.55e-01       & 3.94e-02 & 1.55e-01 \\
$c$ = 0.6  &                  &                &                  &                &          &          \\
         & MSE for [0, 0.5] & MSE for [0, 1] & MAE for [0, 0.5] & MAE for [0, 1] & Max MSE  & Max MAE  \\
RK       & 2.23e-02         & 8.72e-02       & 1.16e-01         & 2.31e-01       & 8.72e-02 & 2.31e-01 \\
Adams3   & 2.23e-02         & 8.72e-02       & 1.16e-01         & 2.31e-01       & 8.72e-02 & 2.31e-01 \\
Uncon    & 7.34e-06         & 2.14e-05       & 2.26e-03         & 3.76e-03       & 2.14e-05 & 3.76e-03 \\
Semi-con & 1.34e-03         & 5.33e-03       & 2.85e-02         & 5.70e-02       & 5.33e-03 & 5.70e-02 \\
Full-con & 2.23e-02         & 8.72e-02       & 1.16e-01         & 2.31e-01       & 8.72e-02 & 2.31e-01 \\
$c$ = 0.7  &                  &                &                  &                &          &          \\
         & MSE for [0, 0.5] & MSE for [0, 1] & MAE for [0, 0.5] & MAE for [0, 1] & Max MSE  & Max MAE  \\
RK       & 3.02e-02         & 1.17e-01       & 1.36e-01         & 2.68e-01       & 1.17e-01 & 2.68e-01 \\
Adams3   & 3.02e-02         & 1.17e-01       & 1.36e-01         & 2.68e-01       & 1.17e-01 & 2.68e-01 \\
Uncon    & 8.47e-06         & 2.92e-05       & 2.41e-03         & 4.32e-03       & 2.92e-05 & 4.32e-03 \\
Semi-con & 1.82e-03         & 7.24e-03       & 3.33e-02         & 6.64e-02       & 7.24e-03 & 6.64e-02 \\
Full-con & 3.02e-02         & 1.17e-01       & 1.36e-01         & 2.68e-01       & 1.17e-01 & 2.68e-01 \\
$c$ = 0.8  &                  &                &                  &                &          &          \\
         & MSE for [0, 0.5] & MSE for [0, 1] & MAE for [0, 0.5] & MAE for [0, 1] & Max MSE  & Max MAE  \\
RK       & 3.94e-02         & 1.52e-01       & 1.55e-01         & 3.05e-01       & 1.52e-01 & 3.05e-01 \\
Adams3   & 3.94e-02         & 1.52e-01       & 1.55e-01         & 3.05e-01       & 1.52e-01 & 3.05e-01 \\
Uncon    & 1.02e-05         & 3.74e-05       & 2.63e-03         & 4.87e-03       & 3.74e-05 & 4.87e-03 \\
Semi-con & 2.38e-03         & 9.44e-03       & 3.81e-02         & 7.58e-02       & 9.44e-03 & 7.58e-02 \\
Full-con & 3.94e-02         & 1.52e-01       & 1.55e-01         & 3.05e-01       & 1.52e-01 & 3.05e-01 \\
$c$ = 0.9  &                  &                &                  &                &          &          \\
         & MSE for [0, 0.5] & MSE for [0, 1] & MAE for [0, 0.5] & MAE for [0, 1] & Max MSE  & Max MAE  \\
RK       & 4.96e-02         & 1.89e-01       & 1.74e-01         & 3.41e-01       & 1.89e-01 & 3.41e-01 \\
Adams3   & 4.96e-02         & 1.89e-01       & 1.74e-01         & 3.41e-01       & 1.89e-01 & 3.41e-01 \\
Uncon    & 1.23e-05         & 4.66e-05       & 2.87e-03         & 5.41e-03       & 4.66e-05 & 5.41e-03 \\
Semi-con & 3.01e-03         & 1.19e-02       & 4.29e-02         & 8.53e-02       & 1.19e-02 & 8.53e-02 \\
Full-con & 4.96e-02         & 1.89e-01       & 1.74e-01         & 3.41e-01       & 1.89e-01 & 3.41e-01 \\
$c$ = 1.0  &                  &                &                  &                &          &          \\
         & MSE for [0, 0.5] & MSE for [0, 1] & MAE for [0, 0.5] & MAE for [0, 1] & Max MSE  & Max MAE  \\
RK       & 6.10e-02         & 2.30e-01       & 1.93e-01         & 3.77e-01       & 2.30e-01 & 3.77e-01 \\
Adams3   & 6.11e-02         & 2.30e-01       & 1.93e-01         & 3.77e-01       & 2.30e-01 & 3.77e-01 \\
Uncon    & 1.50e-05         & 5.72e-05       & 3.14e-03         & 5.99e-03       & 5.72e-05 & 5.99e-03 \\
Semi-con & 3.72e-03         & 1.47e-02       & 4.76e-02         & 9.46e-02       & 1.47e-02 & 9.46e-02 \\
Full-con & 6.11e-02         & 2.30e-01       & 1.93e-01         & 3.77e-01       & 2.30e-01 & 3.77e-01 
\end{longtblr}

\section{Analysis of phase error for the 1-D wave equation}\label{sec:Phase}
First, we analyze the errors introduced only by the spatial discretization. The conservation form of the wave equation can be written as
\begin{align}
\frac{\partial v}{\partial t}+ \frac{\partial J}{\partial x}=0, \quad 0 \leqslant x \leqslant 1, \; 0 \leqslant t \leqslant 1,\label{wave_eq_conservation} \\ 
v(x, 0)=\sin(4 \pi x), \quad 0 \leqslant x \leqslant 1,
\end{align} 
where the flux $J \equiv cv$. The time derivative is approximated using the spatial derivative as
\begin{equation}
\frac{\partial v}{\partial t}\approx   -\frac{J_{j+1/2}-J_{j-1/2}}{\Delta x}.
\label{time_deriva}
\end{equation}
With $v^{n}_{j+1/2}=(\bar{v}^{n}_{j}+\bar{v}^{n}_{j+1})/2$, the Eq.~(\ref{time_deriva}) is written as
\begin{equation}
\left ( \frac{\partial v}{\partial t}\right ) ^{n}_{j}\approx   -c\frac{\bar{v}^{n}_{j+1}-\bar{v}^{n}_{j-1}}{2\Delta x},
\label{time_deriva_n}
\end{equation}
where $\bar{v}$ is the cell-average value, the superscript represents the coordinates of time, and the subscript represents the coordinates of space. 
The cell-average value $\bar{v}$ can be written as
\begin{equation}
\bar{v}_{j}^{n} =\frac{1}{\Delta x} \int_{x_{j-1/2}}^{x_{j+1/2}} v(x,t)dx,\ t=n\Delta t.
\label{cell_average}
\end{equation}
By expanding $v(x,t)$ on the right hand side of Eq.~(\ref{cell_average}) in a Taylor series about $x_{j}$~\cite{lomax2001fundamentals}, one obtains
\begin{equation}
\bar{v}_{j}^{n}=v_{j}^{n}+\frac{\Delta x^{2}}{24}\left(\frac{\partial^{2} v}{\partial x^{2}}\right)_{j}^{n}+\frac{\Delta x^{4}}{1920}\left(\frac{\partial^{4} v}{\partial x^{4}}\right)_{j}^{n}+O\left(\Delta x^{6}\right).
\label{average2center}
\end{equation}
In particular, the initial condition of the equations studied here is given by a sinusoidal function. The even-order derivatives remain sinusoidal without phase changes. Since our main purpose is to investigate the phase lag introduced within a time step, it is reasonable to replace the cell mean $\bar{v}_{j}^{n}$ by the value at the center of the cell $v_{j}^{n}$ in Eq.~(\ref{time_deriva_n}) in error analysis. Thus the error analysis of the employed {\color{black}second-order} finite-volume method is equivalent to the second-order central difference scheme.

To carry out Fourier analysis of the second-order central difference scheme, we employ a single Fourier component
\begin{equation}
v(x,0)=\text{exp} (ikx)
\label{fourier_ini}
\end{equation}
as the initial value. The exact solution of the 1-D wave equation in Eq.~(\ref{wave_eq_conservation}) with initial condition Eq.~(\ref{fourier_ini}) is 
\begin{equation}
v(x,t)=\text{exp}(-ikct)\text{exp}(ikx)=\text{exp}\left [ik(x-ct)\right]. 
\label{exact_fourier}
\end{equation}
By applying the second-order central difference scheme, we can obtain an approximate solution, 
\begin{equation}
\tilde{v}(x,t)=\text{exp}\left(-c\frac{k_{r}}{\alpha}kt \right) \text{exp}\left[ik\left( x-c\frac{k_{i}}{\alpha}t\right)\right],
\label{approximate}
\end{equation}
where $\alpha=k\Delta x$, and $k_{r}=0, k_{i}=\text{sin}\alpha$.

Comparing the above Eq.~(\ref{exact_fourier}) and Eq.~(\ref{approximate}), it can be observed that the propagation speed of the approximate solution for the central difference is $\frac{c}{\alpha} \text{sin}\alpha$, which is always slower than the exact solution $c$. Assuming the solution from the previous time step is exact, the phase lag within a time step for the second-order central difference scheme is
\begin{equation}
d=\left ( 1-\frac{\text{sin}\alpha}{\alpha}\right ) c\Delta t.
\label{phase_lag}
\end{equation}
Roughly speaking, the coarser grids the larger the error will be.

In the following, the time integration scheme and the spatial discretization scheme employed in the test cases will be considered simultaneously. The approximate solution $\tilde{v} ^{n+1}_{j}$ on the $(n+1)\text{-}th$ time layer can be obtained from the solutions at $(n-2)\text{-}th$, $(n-1)\text{-}th$, and $n\text{-}th$ layer as follows:
\begin{equation}
\tilde{v} ^{n+1}_{j}=-\alpha_{0}v^{n-2}_{j}-\alpha_{1}v^{n-1}_{j}-\alpha_{2}v^{n}_{j}+\Delta t\left [ \beta_{0}\left ( \frac{\partial v}{\partial t}  \right ) ^{n-2}_{j} +\beta_{1}\left ( \frac{\partial v}{\partial t}  \right ) ^{n-1}_{j} +\beta_{2}\left ( \frac{\partial v}{\partial t}  \right ) ^{n}_{j} \right ].
\label{wave_lmm}
\end{equation}
By applying Eq.~(\ref{time_deriva_n}) to Eq.~(\ref{wave_lmm}), we  can turn Eq.~(\ref{wave_lmm}) into a four-level explicit scheme:
\begin{equation}
\tilde{v} ^{n+1}_{j}=\left ( \beta_{0}r v^{n-2}_{j-1} -\alpha_{0} v^{n-2}_{j} -\beta_{0}r v^{n-2}_{j+1}\right ) +\left ( \beta_{1}r v^{n-1}_{j-1} -\alpha_{1} v^{n-1}_{j} -\beta_{1}r v^{n-1}_{j+1}\right )+\left ( \beta_{2}r v^{n}_{j-1} -\alpha_{2} v^{n}_{j} -\beta_{2}r v^{n}_{j+1}\right ),
\label{4-layers}
\end{equation}
where $r=\frac{c\Delta t}{2\Delta x}$.

We assume that the solutions from the previous three time steps on the right hand side of Eq.~(\ref{4-layers}) are exact, and aim to formulate the phase error of the numerical solution given by obtained by Eq.~(\ref{4-layers}) at the $(n+1)\text{-}th$ time layer. With the exact solutions given in the following form,  
\begin{equation}
v^{n+\xi }_{j+\ell  }=\text{sin}\left \{ 4\pi\left [ x_{j}+\ell \Delta x-c\left ( t_{n}+\xi\Delta t \right )  \right ]  \right \} ,\quad  \ell =0,\pm 1, \, \xi=-2,-1,0,
\label{v_jn}
\end{equation} 
the right-hand side of Eq.~(\ref{4-layers}) becomes a superposition of 9 sinusoidal functions with the same frequency, which can be denoted as
\begin{equation}
\tilde{v} ^{n+1}_{j}=\sum_{i=1}^{9}a_{i}\text{sin}\left ( 4\pi x+\varphi _{i} \right ),  
\label{9sin}
\end{equation}
where
\begin{gather}
\begin{bmatrix}
  a_{7}& a_{8} &a_{9} \\[5pt]
  a_{4}&a_{5}  &a_{6} \\[5pt]
  a_{1}&a_{2}  &a_{3}
\end{bmatrix}
=\begin{bmatrix}
\beta_{2}r &  -\alpha_{2}&- \beta_{2}r\\[5pt]
\beta_{1}r &  -\alpha_{1}&- \beta_{1}r\\[5pt]
  \beta_{0}r &  -\alpha_{0}&- \beta_{0}r
\end{bmatrix}, \\
\begin{bmatrix}
  \varphi_{7}& \varphi_{8} &\varphi_{9} \\[5pt]
  \varphi_{4}&\varphi_{5}  &\varphi_{6} \\[5pt]
  \varphi_{1}&\varphi_{2}  &\varphi_{3}
\end{bmatrix}
=\begin{bmatrix}
 4\pi\left (-\Delta x-ct\right ) &  4\pi\left (-ct\right ) & 4\pi\left (\Delta x-ct\right )\\[5pt]
 4\pi\left [-\Delta x-c(t-\Delta t)\right ] &  4\pi\left [-c(t-\Delta t)\right ] & 4\pi\left [\Delta x-c(t-\Delta t)\right ]\\[5pt]
 4\pi\left [-\Delta x-c(t-2\Delta t)\right ] &  4\pi\left [-c(t-2\Delta t)\right ] & 4\pi\left [\Delta x-c(t-2\Delta t)\right ]
\end{bmatrix}.
\end{gather}
By using the auxiliary angle formula, the right hand side of the equation can be merged into a sinusoidal function:
\begin{equation}
\sum_{i=1}^{9}a_{i}\text{sin}\left ( 4\pi x+\varphi _{i} \right )=a\text{sin}\left ( 4\pi x+\varphi \right )=a\text{sin}\left [ 4\pi (x+\psi )\right ], 
\label{1sin}
\end{equation}
where 
\begin{equation}
a=\sqrt{\left (   \sum_{i=1}^{9}a_{i}\text{sin}\varphi _{i}\right )^2+\left (   \sum_{i=1}^{9}a_{i}\text{cos}\varphi _{i}\right )^2},
\label{a_witht}
\end{equation}
and $\varphi$ satisfies 
\begin{equation}
\text{sin}\varphi=\frac{\sum_{i=1}^{9}a_{i}\text{sin}\varphi _{i}}{a},\;\text{cos}\varphi=\frac{\sum_{i=1}^{9}a_{i}\text{cos}\varphi _{i}}{a},
\label{varphi}
\end{equation}
and $\psi =\varphi /4\pi$.

As shown in Eqs.~(\ref{a_witht}, \ref{varphi}), the expressions for $a$ and $\varphi$ are composed of the sum of sine and cosine functions with the same frequency. 
In the following, their expressions are simplified with the use of the auxiliary angle formula. Let $\varphi_{i}=-4\pi ct+\Gamma _{i}$, in which $\Gamma _{i}$ is independent of $x$ and $t$ given as follows:
\begin{equation}
\begin{bmatrix}
  \Gamma _{7}& \Gamma _{8} &\Gamma _{9} \\[5pt]
  \Gamma _{4}&\Gamma _{5}  &\Gamma _{6} \\[5pt]
  \Gamma _{1}&\Gamma _{2}  &\Gamma _{3}
\end{bmatrix}
=\begin{bmatrix}
 -4\pi\Delta x &  0 & 4\pi \Delta x\\[5pt]
 4\pi\left (-\Delta x+c\Delta t\right ) &  4\pi\left (c\Delta t\right ) & 4\pi\left (\Delta x+c\Delta t\right )\\[5pt]
 4\pi\left (-\Delta x+2c\Delta t\right ) &  4\pi\left (2c\Delta t\right ) & 4\pi\left (\Delta x+2c\Delta t\right )
\end{bmatrix}.    
\end{equation}
After some straightforward computations, we then have 
\begin{equation}
\sum_{i=1}^{9}a_{i}\text{sin}\varphi _{i}= \sum_{i=1}^{9}a_{i}\text{sin}(-4\pi ct+\Gamma  _{i})=b_{1}\text{sin}(-4\pi ct+\Gamma),  
\label{Gamma_b}
\end{equation}
where
\begin{equation}
b_{1}=\sqrt{\left (\sum_{i=1}^{9}a_{i}\text{sin}\Gamma  _{i}\right )^2+\left (\sum_{i=1}^{9}a_{i}\text{cos}\Gamma _{i}\right )^2},
\label{b_withoutt}
\end{equation}
and $\Gamma $ satisfies that
\begin{equation}
\text{sin}\Gamma =\frac{\sum_{i=1}^{9}a_{i}\text{sin}\Gamma  _{i}}{b_{1}},\; \text{cos}\Gamma =\frac{\sum_{i=1}^{9}a_{i}\text{cos}\Gamma _{i}}{b_{1}}.
\label{Gamma}
\end{equation}
Thus $b_{1}$ and $\Gamma$ are constants for given $\Delta x, \Delta t$ and $c$ for conventional time integration schemes with fixed coefficients. Similarly, we have
\begin{equation}
\sum_{i=1}^{9}a_{i}\text{cos}\varphi _{i}= \sum_{i=1}^{9}a_{i}\text{sin}(\frac{\pi}{2}+4\pi ct-\Gamma  _{i})=b_{2}\text{sin}(4\pi ct+\theta ),  
\label{theta_c}
\end{equation}
where
\begin{equation}
b_{2}=\sqrt{\left [\sum_{i=1}^{9}a_{i}\text{sin}\left ( \frac{\pi}{2}-  \Gamma  _{i}\right ) \right ]^2+\left [\sum_{i=1}^{9}a_{i}\text{cos}\left ( \frac{\pi}{2}-  \Gamma  _{i}\right )\right ]^2}=\sqrt{\left (\sum_{i=1}^{9}a_{i}\text{cos}\Gamma  _{i} \right )^2+\left (\sum_{i=1}^{9}a_{i}\text{sin}\Gamma  _{i}\right )^2}=b_{1},
\label{c=b}
\end{equation}
and $\theta $ satisfies that
\begin{equation}
\text{sin}\theta =\frac{\sum_{i=1}^{9}a_{i}\text{cos}\Gamma  _{i}}{b_{2}},\; \text{cos}\theta =\frac{\sum_{i=1}^{9}a_{i}\text{sin}\Gamma _{i}}{b_{2}},
\label{theta}
\end{equation}
so $\theta$ is a constant as well.

From Eqs.~(\ref{Gamma}, \ref{theta}), it is obvious that 
\begin{equation}
\text{sin}(\theta+\Gamma)=\text{sin}\theta\text{cos}\Gamma+\text{sin}\Gamma\text{cos}\theta=1. 
\end{equation}
As a consequence, $\theta$ and $\Gamma$ satisfy the relationship, $\theta+\Gamma=\pi/2+2k_{1}\pi$, where $k_{1}$ is any integer. After plugging Eqs.~(\ref{Gamma_b}, \ref{theta_c}) into Eq.~(\ref{a_witht}), it follows that $a=b_{1}$, which is independent of $x$ and $t$. That is to say, for the initial condition in Eq.~(\ref{wave_eq_conservation}), if $\Delta x$, $\Delta t$, and $c$ are given, the wave amplitude is a constant during one time-step advance. The value of $a$ for all three learned schemes is approximately 1, indicating no numerical dissipation of the schemes as shown in the results section.

By substituting Eq.~(\ref{Gamma_b}, \ref{theta_c}) into Eq.~(\ref{varphi}), we obtain the following expressions for $\varphi$,
\begin{equation}
\text{sin}\varphi=\text{sin}\left ( -4\pi ct+\Gamma \right ),\;\text{cos}\varphi=\text{cos}\left ( -4\pi ct+\Gamma \right ).
\label{sin_varphi_Gamma}
\end{equation}
From the periodicity of trigonometric functions, we have  
\begin{equation}
\varphi=-4\pi ct+\Gamma+2k_{2}\pi,
\label{varphi_Gamma}
\end{equation}
where $k_{2}$ is any integer. Taking $\psi =\varphi /4\pi$ into consideration, by combining  Eq.~(\ref{9sin}) and  Eq.~(\ref{1sin}), it can be concluded that
\begin{equation}
\tilde{v} ^{n+1}_{j} =a\text{sin}\left[ 4\pi (x_{j}-ct_{n}+\frac{\Gamma}{4\pi})+2k_{2}\pi  \right].  
\label{numerical_sin_n+1}
\end{equation}
It is noted that we aim to examine the phase displacement in one time step for different time integration methods. 
Here $k_{2}=0$ because of small $\Delta t$ and suitable $c$. According to the exact solution (i.e., $v^{n+1}_{j}=\text{sin}\left[4\pi(x_{j}-ct_{n}-c\Delta t) \right]$), the exact phase displacement per one time step is $c\Delta t$.
On the other hand, the phase displacement from the numerical simulation is $-\Gamma /4\pi$.
With the coefficients of the linear multistep method, the values of $-\Gamma /4\pi$ are computed for various time integration methods and plotted in Figure \ref{fig:mean_phase}.
As the coefficients of the learned linear multistep method vary with time, the time-averaged values of $-\Gamma /4\pi$ are plotted.

\section{Detailed results of the 1-D Burgers' equation}\label{sec:Burger' Details }
In this appendix, the errors from different sets of the 1-D Burgers' equation cases are summarized in tables~\ref{tab1} to~\ref{tab5}. The domain is $t\in \left [ 0,100 \right ]$, $x\in \left [ 0,2\pi \right ] $. Tables ~\ref{tab1} to ~\ref{tab3} summarize the error from the learned unconstrained, semi-constrained, and fully-constrained models. Tables \ref{tab4} and \ref{tab5} show the error from the learned unconstrained and semi-constrained models with constant coefficients, which are the approximation of the reserved decimals of the coefficients computed by the corresponding model. 
\begin{table}[ht]
\caption{
Error of the learned semi-constrained model and the Runge-Kutta method for the 20 cases simulated in this work. The first four columns are the mean square error (MSE) and mean absolute error (MAE) from the two methods. The last two columns are the percentages of the MSE/MAE reduction, which is defined as $100 \times \frac{e_{RK}-e_{data-driven}}{e_{RK}}\%$, where $e$ is for MSE or MAE. At the bottom of the table, the p-values from two paired t-tests for evaluating the difference of paired observations in MSE and MAE are presented.
}\label{tab1}
\begin{tabular}{cccccccc}
\toprule
  & Semi-Con MSE & RK MSE & Semi-Con MAE & RK MAE & MSE $\downarrow$ (\%) & MAE $\downarrow$ (\%)\\ \midrule
sample 0     & 0.00578818  & 0.00938566 & 0.0381855   & 0.0446107 & 38.33\%  & 14.40\%  \\
sample 1     & 0.00314997  & 0.00424672 & 0.0314622   & 0.0349126 & 25.83\%  & 9.88\%   \\
sample 2     & 0.0050351   & 0.0068048  & 0.0361549   & 0.0413616 & 26.01\%  & 12.59\%  \\
sample 3     & 0.0133611   & 0.0149293  & 0.0593322   & 0.0616454 & 10.50\%  & 3.75\%   \\
sample 4     & 0.00337585  & 0.00406657 & 0.0357434   & 0.0379615 & 16.99\%  & 5.84\%   \\
sample 5     & 0.0184833   & 0.0201409  & 0.0688312   & 0.0706579 & 8.23\%   & 2.59\%   \\
sample 6     & 0.0216622   & 0.0279999  & 0.0679485   & 0.0753268 & 22.63\%  & 9.80\%   \\
sample 7     & 0.0126621   & 0.0161985  & 0.0586569   & 0.0661358 & 21.83\%  & 11.31\%  \\
sample 8     & 0.00853427  & 0.0100198  & 0.045956    & 0.0476583 & 14.83\%  & 3.57\%   \\
sample 9     & 0.0126793   & 0.0131619  & 0.0558181   & 0.0557141 & 3.67\%   & -0.19\%  \\
sample 10    & 0.00828612  & 0.00828179 & 0.0404443   & 0.04256   & -0.05\%  & 4.97\%   \\
sample 11    & 0.00889011  & 0.0103684  & 0.0441687   & 0.0472441 & 14.26\%  & 6.51\%   \\
sample 12    & 0.00918089  & 0.0122203  & 0.0473384   & 0.052166  & 24.87\%  & 9.25\%   \\
sample 13    & 0.0134826   & 0.0154811  & 0.0523021   & 0.0566761 & 12.91\%  & 7.72\%   \\
sample 14    & 0.00920021  & 0.0130175  & 0.0475373   & 0.0552199 & 29.32\%  & 13.91\%  \\
sample 15    & 0.0127659   & 0.0160335  & 0.0526746   & 0.0572706 & 20.38\%  & 8.03\%   \\
sample 16    & 0.00751947  & 0.00624123 & 0.0430819   & 0.0388038 & -20.48\% & -11.03\% \\
sample 17    & 0.00795443  & 0.00753813 & 0.0464163   & 0.0459236 & -5.52\%  & -1.07\%  \\
sample 18    & 0.00730773  & 0.00725838 & 0.0465312   & 0.0459355 & -0.68\%  & -1.30\%  \\
sample 19    & 0.0115021   & 0.0163913  & 0.0510611   & 0.0586433 & 29.83\%  & 12.93\%  \\ \midrule
\multicolumn{1}{c}{Mean} &
  0.010041047 &
  0.011989284 &
  0.04848224 &
  0.05182138 &
  14.68\% &
  6.17\% \\
\multicolumn{1}{c}{p-value} & 0.000201141 &            & 0.000180005 &           &              &        
\end{tabular}
\end{table}

\begin{table}[ht]
\caption{Error of the learned unconstrained model and the Runge-Kutta method for the 20 cases simulated in this work. The explanations of the headers are the same as Table \ref{tab1}.
}\label{tab2}
\begin{tabular}{cccccccc}
\toprule
 & Un-Con MSE & RK MSE & Un-Con MAE & RK MAE & MSE $\downarrow$ (\%) & MAE $\downarrow$ (\%)\\ \midrule
sample 0  & 0.00601873 & 0.00938566 & 0.0388269 & 0.0446107 & 35.87\%  & 12.97\%  \\
sample 1  & 0.00312457 & 0.00424672 & 0.0313471 & 0.0349126 & 26.42\%  & 10.21\%  \\
sample 2  & 0.00502773 & 0.0068048  & 0.0362164 & 0.0413616 & 26.12\%  & 12.44\%  \\
sample 3  & 0.0128834  & 0.0149293  & 0.0585186 & 0.0616454 & 13.70\%  & 5.07\%   \\
sample 4  & 0.00344257 & 0.00406657 & 0.0359626 & 0.0379615 & 15.34\%  & 5.27\%   \\
sample 5  & 0.0179093  & 0.0201409  & 0.0678392 & 0.0706579 & 11.08\%  & 3.99\%   \\
sample 6  & 0.0206669  & 0.0279999  & 0.0666606 & 0.0753268 & 26.19\%  & 11.50\%  \\
sample 7  & 0.0128389  & 0.0161985  & 0.0589232 & 0.0661358 & 20.74\%  & 10.91\%  \\
sample 8  & 0.00838317 & 0.0100198  & 0.0457445 & 0.0476583 & 16.33\%  & 4.02\%   \\
sample 9  & 0.0125037  & 0.0131619  & 0.0557416 & 0.0557141 & 5.00\%   & -0.04\%  \\
sample 10 & 0.00807299 & 0.00828179 & 0.0396088 & 0.04256   & 2.52\%   & 6.93\%   \\
sample 11 & 0.00899968 & 0.0103684  & 0.0444989 & 0.0472441 & 13.20\%  & 5.81\%   \\
sample 12 & 0.00886234 & 0.0122203  & 0.0469216 & 0.052166  & 27.48\%  & 10.05\%  \\
sample 13 & 0.0128181  & 0.0154811  & 0.0514762 & 0.0566761 & 17.20\%  & 9.17\%   \\
sample 14 & 0.00854321 & 0.0130175  & 0.0464543 & 0.0552199 & 34.37\%  & 15.87\%  \\
sample 15 & 0.012454   & 0.0160335  & 0.0522895 & 0.0572706 & 22.33\%  & 8.70\%   \\
sample 16 & 0.00737186 & 0.00624123 & 0.0428258 & 0.0388038 & -18.12\% & -10.37\% \\
sample 17 & 0.00791533 & 0.00753813 & 0.0463025 & 0.0459236 & -5.00\%  & -0.83\%  \\
sample 18 & 0.007347   & 0.00725838 & 0.0468898 & 0.0459355 & -1.22\%  & -2.08\%  \\
sample 19 & 0.0114697  & 0.0163913  & 0.0508795 & 0.0586433 & 30.03\%  & 13.24\%    \\ \midrule
\multicolumn{1}{c}{Mean}  & 0.009832659 & 0.011989284 & 0.04819638   & 0.05182138 & 15.98\%   & 6.64\%    \\

\multicolumn{1}{c}{p-value} & 0.000147995 &             & 0.000111214 &            &           &   
\end{tabular}
\end{table}

\begin{table}[ht]
\caption{Error of the learned fully-constrained model and the Runge-Kutta method for the 20 cases simulated in this work. The explanations of the headers are the same as Table \ref{tab1}.}\label{tab3}
\begin{tabular}{cccccccc}
\toprule
 & Full-Con MSE & RK MSE & Full-Con MAE & RK MAE & MSE $\downarrow$ (\%) & MAE $\downarrow$ (\%) \\ \midrule
sample 0  & 0.00935032   & 0.00938566  & 0.0445478    & 0.0446107  & 0.38\%        & 0.14\%        \\
sample 1  & 0.00424713   & 0.00424672  & 0.0349076    & 0.0349126  & -0.01\%       & 0.01\%        \\
sample 2  & 0.00679332   & 0.0068048   & 0.0413455    & 0.0413616  & 0.17\%        & 0.04\%        \\
sample 3  & 0.0149495    & 0.0149293   & 0.0616746    & 0.0616454  & -0.14\%       & -0.05\%       \\
sample 4  & 0.00405327   & 0.00406657  & 0.037912     & 0.0379615  & 0.33\%        & 0.13\%        \\
sample 5  & 0.0200868    & 0.0201409   & 0.0705864    & 0.0706579  & 0.27\%        & 0.10\%        \\
sample 6  & 0.0280089    & 0.0279999   & 0.075336     & 0.0753268  & -0.03\%       & -0.01\%       \\
sample 7  & 0.0161669    & 0.0161985   & 0.066082     & 0.0661358  & 0.20\%        & 0.08\%        \\
sample 8  & 0.00999278   & 0.0100198   & 0.0476166    & 0.0476583  & 0.27\%        & 0.09\%        \\
sample 9  & 0.0131692    & 0.0131619   & 0.0557087    & 0.0557141  & -0.06\%       & 0.01\%        \\
sample 10 & 0.00829598   & 0.00828179  & 0.042588     & 0.04256    & -0.17\%       & -0.07\%       \\
sample 11 & 0.0103528    & 0.0103684   & 0.047219     & 0.0472441  & 0.15\%        & 0.05\%        \\
sample 12 & 0.0122202    & 0.0122203   & 0.0521581    & 0.052166   & 0.00\%        & 0.02\%        \\
sample 13 & 0.0154733    & 0.0154811   & 0.0566515    & 0.0566761  & 0.05\%        & 0.04\%        \\
sample 14 & 0.0130123    & 0.0130175   & 0.0551804    & 0.0552199  & 0.04\%        & 0.07\%        \\
sample 15 & 0.0160084    & 0.0160335   & 0.0572202    & 0.0572706  & 0.16\%        & 0.09\%        \\
sample 16 & 0.00622684   & 0.00624123  & 0.0387584    & 0.0388038  & 0.23\%        & 0.12\%        \\
sample 17 & 0.00754246   & 0.00753813  & 0.0459381    & 0.0459236  & -0.06\%       & -0.03\%       \\
sample 18 & 0.00726785   & 0.00725838  & 0.0459242    & 0.0459355  & -0.13\%       & 0.02\%        \\
sample 19 & 0.0164094    & 0.0163913   & 0.0586656    & 0.0586433  & -0.11\%       & -0.04\%       \\ \midrule
\multicolumn{1}{c}{Mean}& 
0.011981383  & 
0.011989284 & 
0.051801035  & 
0.05182138 & 
0.08\%        
& 0.04\% \\

\multicolumn{1}{c}{p-value}& 0.087204851  &             & 0.008416531  &           &              &      
\end{tabular}
\end{table}

\begin{table}[ht]
\caption{Error of the learned unconstrained model with constant coefficients and the Runge-Kutta method for the 20 cases simulated in this work. The explanations of the headers are the same as Table \ref{tab1}.}\label{tab4}
\begin{tabular}{cccccccc}
\toprule
 &Un-Con MSE &RK MSE &Un-Con MAE &RK MAE &MSE $\downarrow$ (\%) &MAE $\downarrow$ (\%) \\ \midrule
sample 0  & 0.00600539 & 0.00938566 & 0.0388007 & 0.0446107 & 36.02\%  & 13.02\%  \\
sample 1  & 0.00312973 & 0.00424672 & 0.0313207 & 0.0349126 & 26.30\%  & 10.29\%  \\
sample 2  & 0.00503258 & 0.0068048  & 0.0362119 & 0.0413616 & 26.04\%  & 12.45\%  \\
sample 3  & 0.0128584  & 0.0149293  & 0.0584661 & 0.0616454 & 13.87\%  & 5.16\%   \\
sample 4  & 0.00344422 & 0.00406657 & 0.0359848 & 0.0379615 & 15.30\%  & 5.21\%   \\
sample 5  & 0.0179302  & 0.0201409  & 0.067894  & 0.0706579 & 10.98\%  & 3.91\%   \\
sample 6  & 0.0207225  & 0.0279999  & 0.0666793 & 0.0753268 & 25.99\%  & 11.48\%  \\
sample 7  & 0.0128853  & 0.0161985  & 0.0590708 & 0.0661358 & 20.45\%  & 10.68\%  \\
sample 8  & 0.00841655 & 0.0100198  & 0.0458378 & 0.0476583 & 16.00\%  & 3.82\%   \\
sample 9  & 0.0125075  & 0.0131619  & 0.0556691 & 0.0557141 & 4.97\%   & 0.08\%   \\
sample 10 & 0.00807902 & 0.00828179 & 0.0395521 & 0.04256   & 2.45\%   & 7.07\%   \\
sample 11 & 0.00895781 & 0.0103684  & 0.0444114 & 0.0472441 & 13.60\%  & 6.00\%   \\
sample 12 & 0.00891361 & 0.0122203  & 0.0469569 & 0.052166  & 27.06\%  & 9.99\%   \\
sample 13 & 0.0129384  & 0.0154811  & 0.051631  & 0.0566761 & 16.42\%  & 8.90\%   \\
sample 14 & 0.00858826 & 0.0130175  & 0.0466744 & 0.0552199 & 34.03\%  & 15.48\%  \\
sample 15 & 0.0125367  & 0.0160335  & 0.0524282 & 0.0572706 & 21.81\%  & 8.46\%   \\
sample 16 & 0.00747684 & 0.00624123 & 0.0430216 & 0.0388038 & -19.80\% & -10.87\% \\
sample 17 & 0.0079493  & 0.00753813 & 0.046369  & 0.0459236 & -5.45\%  & -0.97\%  \\
sample 18 & 0.00733239 & 0.00725838 & 0.0468751 & 0.0459355 & -1.02\%  & -2.05\%  \\
sample 19 & 0.0115947  & 0.0163913  & 0.051101  & 0.0586433 & 29.26\%  & 12.86\%   \\ \midrule
\multicolumn{1}{c}{Mean} &
  0.00986497 &
  0.011989284 &
  0.048247795 &
  0.05182138 &
  15.71\% &
  6.55\% \\
\multicolumn{1}{c}{p-value} & 0.00016039 &            & 0.000120366 &           &              &    
\end{tabular}
\end{table}

\begin{table}[ht]
\caption{
Error of the learned semi-constrained model with constant coefficients and the Runge-Kutta method for the 20 cases simulated in this work. The explanations of the headers are the same as Table \ref{tab1}.}\label{tab5}
\begin{tabular}{cccccccc}
\toprule
 &Semi-Con MSE &RK MSE &Semi-Con MAE &RK MAE &MSE $\downarrow$ (\%) &MAE $\downarrow$ (\%) \\ \midrule
sample 0  & 0.00566226  & 0.00938566  & 0.0385239   & 0.0446107  & 39.67\%   & 13.64\%    \\
sample 1  & 0.00298576  & 0.00424672  & 0.0307164   & 0.0349126  & 29.69\%   & 12.02\%    \\
sample 2  & 0.00484421  & 0.0068048   & 0.0360116   & 0.0413616  & 28.81\%   & 12.93\%    \\
sample 3  & 0.0119549   & 0.0149293   & 0.0569651   & 0.0616454  & 19.92\%   & 7.59\%     \\
sample 4  & 0.0034559   & 0.00406657  & 0.0360134   & 0.0379615  & 15.02\%   & 5.13\%     \\
sample 5  & 0.0171944   & 0.0201409   & 0.0666736   & 0.0706579  & 14.63\%   & 5.64\%     \\
sample 6  & 0.0186734   & 0.0279999   & 0.0641202   & 0.0753268  & 33.31\%   & 14.88\%    \\
sample 7  & 0.0123278   & 0.0161985   & 0.0575442   & 0.0661358  & 23.90\%   & 12.99\%    \\
sample 8  & 0.00813479  & 0.0100198   & 0.0456054   & 0.0476583  & 18.81\%   & 4.31\%     \\
sample 9  & 0.0123754   & 0.0131619   & 0.0559493   & 0.0557141  & 5.98\%    & -0.42\%    \\
sample 10 & 0.00852163  & 0.00828179  & 0.0404007   & 0.04256    & -2.90\%   & 5.07\%     \\
sample 11 & 0.00895405  & 0.0103684   & 0.0446266   & 0.0472441  & 13.64\%   & 5.54\%     \\
sample 12 & 0.00823713  & 0.0122203   & 0.0459845   & 0.052166   & 32.59\%   & 11.85\%    \\
sample 13 & 0.0119645   & 0.0154811   & 0.0501744   & 0.0566761  & 22.72\%   & 11.47\%    \\
sample 14 & 0.00753242  & 0.0130175   & 0.0444447   & 0.0552199  & 42.14\%   & 19.51\%    \\
sample 15 & 0.0113857   & 0.0160335   & 0.0506332   & 0.0572706  & 28.99\%   & 11.59\%    \\
sample 16 & 0.00790158  & 0.00624123  & 0.0443181   & 0.0388038  & -26.60\%  & -14.21\%   \\
sample 17 & 0.00815222  & 0.00753813  & 0.0465709   & 0.0459236  & -8.15\%   & -1.41\%    \\
sample 18 & 0.00773376  & 0.00725838  & 0.0481125   & 0.0459355  & -6.55\%   & -4.74\%    \\
sample 19 & 0.0104727   & 0.0163913   & 0.0491817   & 0.0586433   & 36.11\%   & 16.13\%    \\ \midrule
\multicolumn{1}{c}{Mean} &
  0.009423226 &
  0.011989284 &
  0.04762852 &
  0.05182138 &
  18.09\% &
  7.48\% \\
\multicolumn{1}{c}{p-value} & 0.000341949 &            & 0.000350297 &           &              &       
\end{tabular}
\end{table}

\FloatBarrier
\end{appendices}
\end{document}